\documentclass[12pt,reqno]{amsart}

\usepackage[utf8]{inputenc}
\usepackage{amsfonts}
\usepackage{amsmath}
\usepackage{amssymb}
\usepackage{amsthm}
\usepackage{mathrsfs}
\usepackage{color}
\usepackage[english]{babel}
\usepackage[T1]{fontenc}
\usepackage{url}
\usepackage{graphicx}
\usepackage[pdftex,colorlinks,backref=page,citecolor=blue]{hyperref}
\usepackage{caption}
\usepackage{enumerate}
\usepackage{hyphenat}
\usepackage{float}
\usepackage{indentfirst}
\usepackage[export]{adjustbox}
\usepackage{tikz}
\usepackage[all]{xy}
\usetikzlibrary{matrix,arrows,decorations.pathmorphing}
\usepackage{booktabs}
\usepackage{array}
\usepackage{csquotes}
\usepackage{todonotes}
\usepackage{comment}
\usepackage{multirow}

\usepackage{float}

\usepackage{vcell}

\newcolumntype{P}[1]{>{\centering\arraybackslash}p{#1}}

\allowdisplaybreaks

\newtheorem{theorem}{Theorem}[section]
\newtheorem{lemma}[theorem]{Lemma}
\newtheorem{proposition}[theorem]{Proposition}
\newtheorem{corollary}[theorem]{Corollary}
\newtheorem{lemma-definition}[theorem]{Lemma--Definition}

\theoremstyle{definition}
\newtheorem{remark}[theorem]{Remark}

\newtheorem{definition}[theorem]{Definition}

\def\PP{{\mathbb P}}

\makeatletter
\def\l@subsection{\@tocline{2}{0pt}{2.5pc}{5pc}{}}
\makeatother

\author{Yongnam Lee}
\address{Center for Complex Geometry, Institute for Basic Science (IBS), 55 Expo-ro, Yuseong-gu, Daejeon, 34126 Korea}
\email{ynlee@ibs.re.kr}

\author{Elisa Scattone}
\address{Dipartimento di Matematica, Universit\`a degli studi di Roma Tor Vergata, Via della Ricerca Scientifica 1, 00133, Roma, Italy}
\email{scattone@mat.uniroma2.it}

\author{Luca Schaffler}
\address{Dipartimento di Matematica e Fisica, Universit\`a degli Studi Roma Tre, Largo San Leonardo Murialdo 1, 00146, Roma, Italy}
\email{luca.schaffler@uniroma3.it}

\title{K3 tails in the KSBA moduli of Kunev surfaces}


\begin{document}

\begin{abstract}
In the context of compact moduli spaces of Koll\'ar, Shepherd-Barron, and Alexeev stable surfaces, a degeneration can acquire an exceptional K3 component in the limit, or more simply a \emph{K3 tail}. We study this phenomenon for Kunev surfaces: minimal surfaces of general type with $K^2=p_g=1$ that arise as $(\mathbb{Z}/2\mathbb{Z})^2$-covers of the projective plane branched along two general cubics and a line. We classify the K3 tails that occur when the branch curve degenerates to a singular curve with a non-log canonical isolated semi-quasihomogeneous singularity, and describe the corresponding boundary strata. This builds on the Kerr--Laza classification of K3 tails arising from isolated semi-quasihomogeneous hypersurface singularities.
\end{abstract}

\maketitle

\section{Introduction}

In \cite{Kun77}, Kunev constructed certain minimal surfaces of general type satisfying $p_g=K^2=1$ that gave a counterexample to the infinitesimal Torelli theorem, answering a question of Griffiths \cite[Problem~7.1]{Gri70}. More in general, Catanese \cite{Cat79} and independently Todorov \cite{Tod80} gave a description of all minimal surfaces of general type with $K^2=p_g=1$. The canonical model of such a surface $S$ is a weighted complete intersection of bidegree $(6,6)$ in the weighted projective space $\mathbb{P}(1,2,2,3,3)$, which form an $18$-dimensional moduli space. In particular, $S$ can be viewed as degree $4$ cover of $\mathbb{P}^2$ via the restriction of the projection $\mathbb{P}(1,2,2,3,3)\dashrightarrow\mathbb{P}(1,2,2)\cong\mathbb{P}^2$. The case where this cover is Galois recovers Kunev's example: that is, $S$ is a $(\mathbb{Z}/2\mathbb{Z})^2$-cover of $\mathbb{P}^2$ branched along two smooth cubics $C_1,C_2$ and a line $L$ such that $C_1+C_2+L$ has at worst nodal singularities (a Galois $\mathbb{Z}/4\mathbb{Z}$-cover cannot occur by \cite[Theorem~3]{Cat79}).

In particular, Kunev's surfaces form a $12$-dimensional moduli space $\mathrm{M}_\mathrm{K}$. A partial compactification of $\mathrm{M}_{\mathrm{K}}$ was constructed by Usui in \cite{Usu89} considering degenerations where $C_1+C_2$ acquires ADE singularities and $L$ is any line: see Section~\ref{sec:previous-result} for a brief summary. In this paper we turn our attention to the projective compactification $\overline{\mathrm{M}}_\mathrm{K}$ parametrizing stable surfaces in the sense of Koll\'ar, Shepherd-Barron, and Alexeev \cite{KSB88,Ale96,Kol23}. KSBA-stable limits of Kunev surfaces such that $\left(\mathbb{P}^2,\frac{1}{2}(C_1+C_2+L)\right)$ is log canonical were studied in \cite[Section~5]{FPR17}. In the current paper we go beyond this by analyzing the stable limits when $\left(\mathbb{P}^2,\frac{1}{2}(C_1+C_2+L)\right)$ acquires non-log canonical singularities. (For the Gorenstein stable limits of the more general Catanese--Todorov surfaces see \cite{DR22}.)

Our motivation also sprouts from the case of I-surfaces (smooth minimal general type surfaces with $K^2=1,p_g=2$) as studied in \cite{GPSZ24}. By considering one-parameter degenerations of I-surfaces acquiring an isolated non-log canonical singularity, the stable replacements of such limits describe eight boundary divisors of the KSBA compactification. Such stable limits consist of two surfaces glued along a $\mathbb{P}^1$, and the exceptional surfaces of the weighted blows up yielding the stable limits are K3 surfaces in weighted $3$-dimensional projective space. Our objective is to find analogous stable limits parametrized by $\overline{\mathrm{M}}_{\mathrm{K}}$, and obtain a complete classification of them.

We start from the isolated semi-quasihomogeneous curve singularities $Z_{12},W_{15},Z_{15},N_{16}$ from \cite{Arn76}. These occur as isolated singularity $\xi=[1:0:0]$ of the branch curve $C_1+C_2+L$, or more precisely $\xi\in(C_1+C_2)\setminus L$: see Lemma--Definition~\ref{def:Z_12-W_15-Z_1,0-N_16-arrangements}. These induce a non-log canonical singular point for the $(\mathbb{Z}/2\mathbb{Z})^2$-cover $S_0$. By studying one-parameter degenerations of Kunev surfaces to $S_0$, we identify several loci in the boundary of $\overline{\mathrm{M}}_{\mathrm{K}}$ parametrizing the corresponding stable limits. The following theorem summarizes the results and constructions of Sections~\ref{sec:non-log-canonical-singularities-Z12-W15-Z15-N16}--\ref{sec:stable-non-normal-K3-limits-from-N16}; more detailed references are provided in the discussion that follows.

\begin{theorem}
\label{thm:first-main-result}
There exist irreducible boundary strata
\[
D(Z_{12}),D(W_{15}),D(Z_{15}),D(N_{16})\subseteq\overline{\mathrm{M}}_{\mathrm{K}},
\]
of codimension $1,1,1,4$, respectively. The generic surface parametrized by $D(Z_{12})$, $D(W_{15})$, $D(Z_{15})$ is in the form $\widetilde{Z}\cup\widetilde{Y}$, where $\widetilde{Z}$ is a rational surface with only finite cyclic quotient singularities, $\widetilde{Z}\cap\widetilde{Y}\cong\mathbb{P}^1$, and $\widetilde{Y}$ is an ADE K3 surface isomorphic to the complete intersection $X_{n,m}$ of two hypersurfaces of degrees $n$ and $m$ in weighted projective $4$-dimensional space as follows:
\begin{itemize}

\item For $Z_{12}$, $X_{10,12}\subseteq\mathbb{P}(1,4,6,5,6)$;

\item For $W_{15}$, $X_{6,6}\subseteq\mathbb{P}(1,2,3,3,3)$;

\item For $Z_{15}$, $X_{6,8}\subseteq\mathbb{P}(1,2,4,4,3)$.

\end{itemize}
The generic surface parametrized by $D(N_{16})$ is a non-normal irreducible surface whose normalization is an ADE K3 surface of the form $X_{4,6}\subseteq\mathbb{P}(1,2,2,2,3)$.
\end{theorem}

The codimension $2$ K3 surfaces $X_{n,m}$ in the statement are, respectively, \#42, \#13, \#19, \#8 in Iano-Fletcher's \cite[Table~2]{IF00}. The boundary divisors $D(Z_{12}),D(W_{15}),D(Z_{15})$ appear to be new. The boundary stratum $D(N_{16})$ already appeared in \cite[Example~5.4]{FPR17}, so our only contribution in this case is describing how the surfaces generically parametrized by it arise from the $N_{16}$ singularity.

We now discuss more in detail the construction of such stable limit surfaces. Consider again the non-log canonical surface $S_0$ given by the $(\mathbb{Z}/2\mathbb{Z})^2$-cover of $\mathbb{P}^2$ branched along $C_1+C_2+L$ with a singularity $\Sigma\in\{Z_{12},W_{15},Z_{15},N_{16}\}$ at $\xi\in C_1+C_2$. In Definition~\ref{def:1-parameter-degeneration-Sigma} we construct explicit one-parameter degenerations $\mathcal{S}\rightarrow\Delta$, where $\Delta=\mathrm{Spec}(R)$ and $R=\mathbb{C}[[t]]$, with central fiber $S_0$, generic fiber a smooth Kunev surface, and for which we can compute the stable replacement of the central fiber. These have a fiberwise $(\mathbb{Z}/2\mathbb{Z})^2$-covering structure $\mathcal{S}\rightarrow\Delta\times\mathbb{P}^2=:\mathcal{X}$ branched along $\mathcal{B}\subseteq\mathcal{X}$, so that the central fiber $\left(X_0,\frac{1}{2}B_0\right)$ of the one parameter family $\left(\mathcal{X},\frac{1}{2}\mathcal{B}\right)\rightarrow\Delta$ is $\left(\mathbb{P}^2,\frac{1}{2}(C_1+C_2+L)\right)$. By the theory in \cite{AP12}, the stable replacement of the central fiber of $\mathcal{S}\rightarrow\Delta$ is an appropriate $(\mathbb{Z}/2\mathbb{Z})^2$-cover of the stable replacement for $\left(\mathcal{X},\frac{1}{2}\mathcal{B}\right)\rightarrow\Delta$, so we compute the latter.

In Proposition~\ref{prop:central-fiber-after-weighted-blow-up} we show that an appropriate weighted blow up $\mathcal{X}'\rightarrow\mathcal{X}$ with center $(0,\xi)$ and exceptional divisor $Y$ yields the stable limit for $\Sigma=Z_{12}$. For $\Sigma=W_{15}$, we follow the weighted blow up by the contraction $\mathcal{X}'\rightarrow\mathcal{X}''$ of a curve in the strict transform of $X_0$, and similarly for $Z_{15}$ (see Proposition~\ref{prop:contraction-W15-Z1,0-case}). For $\Sigma=N_{16}$, it is necessary to contract the whole strict transform of $X_0$ (see Proposition~\ref{prop:stable-replacement-for-N16}) to obtain the stable limit. The boundary strata $D(Z_{12}),D(W_{15}),D(Z_{15})$ (resp. $D(N_{16})$) are constructed in Definition~\ref{def:boundary-stratum-associated-to-Z122} (resp. Definition~\ref{def:boundary-stratum-associated-to-N16}), and their dimension count is proved in Theorem~\ref{thm:divisor-for-Z_12^2} (resp. Proposition~\ref{prop:dimension-count-for-D(N16)}). The $(\mathbb{Z}/2\mathbb{Z})^2$-covers $\widetilde{Y}\rightarrow Y$ of the exceptional divisors $Y\subseteq\mathcal{X}'$ are computed in Section~\ref{sec:description-of-the-cover-of-Y}, which is were we identify them as codimension $2$ K3 surfaces in weighted projective $4$-dimensional space.

Finally, we ask whether there are other ways of letting $C_1+C_2+L$ acquire a semi-quasihomogeneous singularity $\Sigma$ such that the exceptional surface $\widetilde{Y}$ in the birational modification $\mathcal{S}'\rightarrow\mathcal{S}$ resolving the non-log canonical singularity is an ADE K3 surface. The appropriate formulation of one-parameter family $\mathcal{S}\rightarrow\Delta$ that generalizes the ones we studied for $\Sigma\in\{Z_{12},W_{15},Z_{15},N_{16}\}$ can be found in Definition~\ref{def:K3-tail-in-moduli-of-Kunev-surfaces}. In this sense, we prove that these are the only cases that yield an ADE K3 surface $\widetilde{Y}$.

\begin{theorem}[Theorem~\ref{thm:classification-of-K3-tails}]
\label{thm:second-main-result}
Let $\mathcal{S}\rightarrow\Delta$ be a one-parameter degeneration of Kunev surfaces as in Definition~\ref{def:K3-tail-in-moduli-of-Kunev-surfaces}, and let $\mathcal{S}'\rightarrow\mathcal{S}$ be the associated blow up of the isolated non-log canonical singularity of the central fiber. Let $\widetilde{Y}$ be the corresponding exceptional divisor, which is an ADE K3 surface. Then, $\widetilde{Y}$ is a codimension $2$ complete intersection in a $4$-dimensional weighted projective space and the only possible cases are $\#8$ ($N_{16}$), $\#13$ ($W_{15}$), $\#19$ ($Z_{15}$), and $\#42$ ($Z_{12}$) in \cite[Table~2]{IF00}.
\end{theorem}

We give an outline of the proof of this classification. We first show that $\xi\in L$ does not give rise to a K3 surface $\widetilde{Y}$ (the cases were the singularity is on $L$ were considered by Usui in \cite{Usu89}). If $\xi\notin L$, then quotienting $\widetilde{Y}$ by an appropriate involution $\iota$ from the $(\mathbb{Z}/2\mathbb{Z})^2$-action gives a K3 surface that is the exceptional surface of the weighted blow up of the non-log canonical semi-quasihomogeneous singularity of the double cover of $\mathbb{P}^2$ branched along $C_1+C_2$. Then, the K3 surface $\widetilde{Y}/\iota$ must fall within the Kerr--Laza classification in \cite{KL25}. Here, we find the four cases in Theorem~\ref{thm:first-main-result}, and we exclude the remaining ones using the constraints imposed by the use of two cubic curves $C_1,C_2$, plus the combinatorics of the list of codimension $2$ ADE K3 surfaces in \cite[Table~2]{IF00}.

\begin{remark}
In complete generality, explicitly describing the boundary of a given KSBA compactification is a difficult problem. For moduli of stable surfaces, one geometric feature that we can aim to understand is the appearance of K3 tails and the structure of the corresponding boundary strata, as we have done in this paper for Kunev surfaces. We hope that the methods developed here can also be applied to study the K3 tails and their associated boundary strata in other KSBA moduli spaces.
\end{remark}


\subsection*{Acknowledgments}
We would like to thank Fabrizio Catanese, Patricio Gallardo, Rita Pardini, and Gregory Pearlstein for helpful conversations. We thank SCGP and ICERM for their hospitality and support during our visits for the conferences ``2nd Simons Math Summer Workshop: Moduli'' and ``Moduli of K3 surfaces''.

\subsection*{Funding}
Y. Lee is supported by the Institute for Basic Science
(IBS-R032-D1). L. Schaffler was supported by the projects ``Programma per Giovani Ricercatori Rita Levi Montalcini'', PRIN2020KKWT53 ``Curves, Ricci flat Varieties and their Interactions'', and PRIN 2022 ``Moduli Spaces and Birational Geometry'' – CUP E53D23005790006. E. Scattone and L. Schaffler are members of the INdAM group GNSAGA and are supported by the INdAM - GNSAGA Project ``Birational geometry of compact moduli spaces via stable pairs and foliations'' - CUP E53C25002010001.

\section{Preliminaries}

\subsection{Kunev surfaces}

Let $C_1,C_2$ be two smooth cubic curves in $\mathbb{P}^2$ and $L$ a line such that $C_1+C_2+L$ only has nodal singularities. By letting $(\mathbb{Z}/2\mathbb{Z})^2=\{e,a,b,c\}$ with $e$ the identity element, let $S\rightarrow\mathbb{P}^2$ be the $(\mathbb{Z}/2\mathbb{Z})^2$-cover with the following building data (see \cite[Section~2]{Par91}):
\[
D_a=C_1,~D_b=C_2,~D_c=L.
\]
Then, $S$ is a complex smooth minimal algebraic surface of general type called \emph{Kunev surface}. Its numerical invariants are $K_S^2=p_g(S)=1$ and $q(S)=0$ (we note that the vanishing of $q(S)$ is automatic from the condition $K_S^2<2p_g(S)$ by \cite{Deb82}). Alternatively, we can view $S$ as a $(6,6)$ codimension $2$ complete intersection in the weighted projective space $\mathbb{P}(1,2,2,3,3)$ as follows. Up to an automorphism of $\mathbb{P}^2$, we can suppose that $L=V(x)$. Then $S\subseteq\mathbb{P}(1,2,2,3,3)$ is given by
\begin{displaymath}
\left\{ \begin{array}{ll}
u^2=F_1(x^2,y,z)\\
v^2=F_2(x^2,y,z),
\end{array} \right.
\end{displaymath}
where $C_1=V(F_1)$, $C_2=V(F_2)$, and the coordinates in $\mathbb{P}(1,2,2,3,3)$ are $[x:y:z:u:v]$. The covering map is the restriction to $S$ of the projection $\mathbb{P}(1,2,2,3,3)\dashrightarrow\mathbb{P}(1,2,2)$, followed by the isomorphism $\mathbb{P}(1,2,2)\cong\mathbb{P}^2$.

\subsection{The moduli space of Kunev surfaces}
\label{sec:mod-space-kunev-surf}

We construct the moduli spaces of Kunev surfaces via their realization as $(\mathbb{Z}/2\mathbb{Z})^2$-cover of $\mathbb{P}^2$ branched along two cubics and a line. Let us consider two copies of $\mathbb{P}^9=\mathbb{P}(H^0(\mathbb{P}^2,\mathcal{O}_{\mathbb{P}^2}(3)))$, one with coordinates $[a_{300}:a_{210}:\ldots:a_{003}]$ and the second one with coordinates $[b_{300}:b_{210}:\ldots:b_{003}]$. Also, let $\mathbb{P}^2=\mathbb{P}(H^0(\mathbb{P}^2,\mathcal{O}_{\mathbb{P}^2}(1)))$ with coordinates $[c_{100}:c_{010}:c_{001}]$. Let $\mathcal{U}\subseteq\mathbb{P}_a^9\times\mathbb{P}_b^9\times\mathbb{P}_c^2$ be the dense open subset such that the plane curve given by the vanishing locus of the following polynomial:
\[
f_{abc}(x,y,z):=\sum_{\substack{0\leq i,j,k\leq3\\i+j+k=3}}a_{ijk}x^iy^jz^k\cdot\sum_{\substack{0\leq i,j,k\leq3\\i+j+k=3}}b_{ijk}x^iy^jz^k\cdot(c_{100}x+c_{010}y+c_{001}z)
\]
consists of the union of two smooth cubics $C_a$, $C_b$ and a line $L_c$ such that $C_a+C_b+L_c$ has only nodal singularities. Let us define the incidence variety
\[
\mathcal{C}:=\{(a,b,c,[x:y:z])\in\mathcal{U}\times\mathbb{P}^2\mid f_{abc}(x,y,z)=0\}.
\]
Then, the fiber of the family $(\mathcal{U}\times\mathbb{P}^2,\mathcal{C})\rightarrow\mathcal{U}$ over $(a,b,c)$ is given by $(\mathbb{P}^2,C_a+C_b+L_c)$. Let $\mathcal{S}\rightarrow\mathcal{U}\times\mathbb{P}^2$ be the appropriate $(\mathbb{Z}/2\mathbb{Z})^2$-cover branched along $\mathcal{C}$, so that the fibers of $\mathcal{S}\rightarrow\mathcal{U}$ recover are all the possible Kunev surfaces. To construct their moduli space, it is well-known that the isomorphism class of a Kunev surface $S$ is uniquely determined by the isomorphism class of the pair $(\mathbb{P}^2,C_a+C_b+L_c)$, where the two cubics can possibly be switched (see for instance \cite{Usu89}). So, let us consider $\mathrm{SL}_3$ acting on $\mathbb{P}^9\times\mathbb{P}^9\times\mathbb{P}^2$ by linear change of coordinates on the polynomial $f_{abc}(x,y,z)$. Also, consider the $S_2$-action given by permuting the two copies of $\mathbb{P}^9$. Then, this induces an action $\mathrm{SL}_3\times S_2\curvearrowright\mathcal{U}$, so that we can take the quotient
\[
\mathrm{M}_{\mathrm{K}}:=\mathcal{U}/(\mathrm{SL}_3\times S_2).
\]
We refer to $\mathrm{M}_{\mathrm{K}}$ as the \emph{moduli space of Kunev surfaces}. It follows from this construction that
\[
\dim(\mathrm{M}_{\mathrm{K}})=\dim(\mathcal{U})-\dim(\mathrm{SL}_3\times S_2)=20-8=12.
\]
We point out that $\mathrm{M}_{\mathrm{K}}$ is parametrizing isomorphism classes of Kunev surfaces $S$ or, equivalently, isomorphism classes of pairs in the form $\left(\mathbb{P}^2,\frac{1}{2}(C_a+C_b+L_c)\right)$, where the two cubics are unmarked. The reason for the choice of the coefficient $\frac{1}{2}$ is because
\[
K_S\sim_\mathbb{Q}\pi^*\left(K_{\mathbb{P}^2}+\frac{1}{2}(C_a+C_b+L_c)\right),
\]
where $\pi\colon S\rightarrow\mathbb{P}^2$ is the $(\mathbb{Z}/2\mathbb{Z})^2$-covering morphism.

\subsection{The KSBA compactification of \texorpdfstring{$\mathrm{M}_{\mathrm{K}}$}{M_K}}
\label{sec:construction-of-the-compactification}

Let $X$ be a $d$-dimensional variety and $\mathbf{b}B:=\sum_ib_iB_i$ a $\mathbb{Q}$-divisor on $X$ with $B_i$ prime divisors and $0<b_i\leq1$. Then, the pair $(X,B)$ is called \emph{stable} provided it has \emph{semi-log canonical singularities} \cite[Definition--Lemma~5.10]{Kol13} and $K_X+B$ is ample. If $B=0$, then $X$ is called a \emph{stable variety}. By \cite[Section~8]{Kol23} we have a moduli stack $\mathcal{SP}(\mathbf{b},d,c)$ parametrizing the isomorphism classes of stable pairs $(X,\mathbf{b}B)$ such that $d=\dim(X)$ and $(K_X+\mathbf{b}B)^d=c$.

In the current paper, we are interested in the case $\mathcal{SP}(0,2,1)$: the family of Kunev surfaces $\mathcal{S}\rightarrow\mathcal{U}$ constructed in Section~\ref{sec:mod-space-kunev-surf} induces then a morphism $\mathcal{U}\rightarrow\mathcal{SP}(0,2,1)$, which factors through $\mathrm{M}_{\mathrm{K}}$. Let $\overline{\mathcal{M}}_{\mathrm{K}}$ be the closure of the image of $\mathcal{U}$ inside $\mathcal{SP}(0,2,1)$ and let $\overline{\mathrm{M}}_{\mathrm{K}}$ be its corresponding projective coarse moduli space. We refer to $\overline{\mathrm{M}}_{\mathrm{K}}$ as the \emph{KSBA compactification of the moduli space of Kunev surfaces}. We note that $\overline{\mathrm{M}}_{\mathrm{K}}$ is parametrizing stable limits of Kunev surfaces, or equivalently stable limits of pairs $\left(\mathbb{P}^2,\frac{1}{2}(C_1+C_2+L)\right)$. This is because the theory of abelian covers is compatible with taking stable limits: see \cite{AP12}.

\subsection{Previous results}
\label{sec:previous-result}

Although we will not use them, for completeness we review some previous results about degenerations of Kunev surfaces. We start from \cite{Usu89}, where Usui described a partial compactification $\mathrm{M}_{\mathrm{K}}^*$ of $\mathrm{M}_{\mathrm{K}}$ that we now briefly recount. For reduced sextic curves in $\PP^2$, curves with at worst simple singularities coincide with GIT stable curves with respect to the action of $\mathrm{SL}_3$ \cite{Shah80}. If we consider the special case where the sextic is the union of two cubics $C_1$ and $C_2$, then the corresponding open subset $U\subseteq\mathrm{Sym}^2|\mathcal{O}_{\mathbb{P}^2}(3)|$ parametrizes ADE K3 surfaces given by the double cover of $\mathbb{P}^2$ branched along $C_1+C_2$. For any such $C_1+C_2$, Usui considered the degenerations of Kunev surfaces that arise as $(\mathbb{Z}/2\mathbb{Z})^2$-covers of $\mathbb{P}^2$ branched along $C_1+C_2+L$ for any choice of line $L\in|\mathcal{O}_{\mathbb{P}^2}(1)|$. In particular, Usui's partial compactification $\mathrm{M}_{\mathrm{K}}^*$ is $(U\times|\mathcal{O}_{\mathbb{P}^2}(1)|)/\mathrm{SL}_3$.

The minimal models of the minimal resolution of surfaces in $\mathrm{M}_{\mathrm{K}}^*$ can be Kunev surfaces, proper elliptic surfaces with $p_g=1, q=0$ with one multiple fiber with multiplicity $2$, proper elliptic surfaces with $p_g=q=1$ without multiple fibers, K3 surfaces, or abelian surfaces. For instance, if the line $L$ passes through one of the nine intersection points of two smooth general cubics, then $S$ has one $\frac{1}{4}(1,1)$-singularity, and the minimal model of the minimal resolution of $S$ is a proper elliptic surfaces with $p_g=1, q=0$ with one multiple fiber of multiplicity $2$. The exceptional $(-4)$-curve of the minimal resolution is a $4$-multisection of this elliptic surface. Such singular surfaces $S$ are stable, and they are hence parametrized by points in $\overline{\mathrm{M}}_{\mathrm{K}}$. The closure of this locus gives a boundary divisor $D\left(\frac{1}{4}(1,1)\right)\subseteq\overline{\mathrm{M}}_{\mathrm{K}}$. We remark that there are also points in Usui's partial compactification $\mathrm{M}_{\mathrm{K}}^*$ parametrizing surfaces whose singularities are not semi-log canonical. For instance, we could let two smooth cubics $C_1,C_2$ form an $A_5$ singularity and let $L$ be the tangent to this singularity so that in the affine patch $x\neq0$, $C_1+C_2+L$ is given by $z(y^6+z^2)=0$. Then, $\left(\mathbb{P}^2,\frac{1}{2}(C_1+C_2+L)\right)$ is not log canonical because its log canonical threshold is $\frac{4}{9}$.

From a different prospective, \cite[Section~5]{FPR17} considered the case where the pair $\left(\mathbb{P}^2,\frac{1}{2}(C_1+C_2+L)\right)$ is log canonical. For instance, the following is \cite[Proposition~5.1]{FPR17}.

\begin{proposition}
Let $S$ be a stable surface parametrized by a point in the KSBA compactification $\overline{\mathrm{M}}_{\mathrm{K}}$. Suppose that we have a $(\mathbb{Z}/2\mathbb{Z})^2$-cover morphism $\pi\colon S\rightarrow\mathbb{P}^2$ branched along $C_1+C_2+L$, where the pair $(\mathbb{P}^2,C_1+C_2+L)$ is not necessarily log smooth. Then, $S$ is $2$-Gorenstein if and only if $C_1\cap C_2\cap L\neq\emptyset$. Furthermore, $S$ is non-normal above an irreducible curve $\Gamma$ of $\mathbb{P}^2$ if and only if $\Gamma$ appears with coefficient $1$ in $\Delta=\frac{1}{2}(C_1+C_2+L)$.
\end{proposition}

As we will see, in the current paper we determine several other geometric boundary strata in $\overline{\mathrm{M}}_{\mathrm{K}}$ different from these ones previously known by computing certain stable limits when $\left(\mathbb{P}^2,\frac{1}{2}(C_1+C_2+L)\right)$ acquires a non-log canonical singularity.

\section{Non lc limits of Kunev surfaces and the singularities \texorpdfstring{$Z_{12}$}{Z 12}, \texorpdfstring{$W_{15}$}{W 15}, \texorpdfstring{$Z_{15}$}{Z 15}, \texorpdfstring{$N_{16}$}{N 16}}
\label{sec:non-log-canonical-singularities-Z12-W15-Z15-N16}

In this section we introduce specific non-log canonical isolated surface singularities which can be acquired by a degeneration of Kunev surfaces.

\subsection{The singularities}
\label{sec:the-singularities-considered}

Our construction is modeled around the following normal forms of plane curve singularities:
\begin{itemize}

\item $Z_{12}\colon yz(y^3+z^2)+ay^3z^2=0$;

\item $W_{15}$ or $W_{1,0}\colon(y^6+ay^3z^2+z^4)+by^4z^2=0$, $a^2\neq4$;

\item $Z_{15}$ or $Z_{1,0}\colon y(y^6+ay^2z^2+z^3)+by^6z=0$, $4a^3+27\neq0$;

\item $N_{16}\colon yz(y^3+ay^2z+byz^2+z^3)+cy^3z^3$, $ab\neq9$, $4(a^3+b^3)+
27-a^2b^2-18ab\neq0$.

\end{itemize}
These normal forms of singularities (up to relabeling of the variables, reordering, and factoring of the monomials) are from \cite{Arn76}. In parentheses we highlighted the quasihomogeneous parts of the singularities for the following choices of weights $(p,q)$ for $(y,z)$ and weighted degree $d_\Sigma$:
\begin{itemize}
\item $Z_{12}$: $(p,q)=(2,3)$, $d_\Sigma=11$;
\item $W_{15}$: $(p,q)=(2,3)$, $d_\Sigma=12$;
\item $Z_{15}$: $(p,q)=(1,2)$, $d_\Sigma=7$;
\item $N_{16}$: $(p,q)=(1,1)$, $d_\Sigma=5$.
\end{itemize}

Writing a local equation of $C_1+C_2$ as the product of the local equations of its branches, the prescribed tangency conditions give the corresponding weighted principal parts.  The asserted analytic types then follow from Arnol'd's classification: see \cite[\S1--2, pp.~93--94 and 102--103, $20_1, 23_2, 30_1$, and $48$]{Arn76}.

\begin{lemma-definition}
\label{def:Z_12-W_15-Z_1,0-N_16-arrangements}
Let $\xi\in\mathbb{P}^2$ be the point $[1:0:0]$. Let $C_1,C_2$ be two cubic curves as follows:
\begin{enumerate}

\item $C_1$ has a node at $\xi$, and $C_2$ has a cusp at $\xi$ whose tangent is one of the two tangent lines of $C_1$. Then, locally analytically $(C_1+C_2,\xi)$ is of type $Z_{12}$.

\item Both $C_1,C_2$ have cusps at $\xi$ with the same tangent line and $I_\xi(C_1,C_2)=6$. Then, locally analytically $(C_1+C_2,\xi)$ is of type $W_{15}$.

\item The curve $C_1$ has a node at $\xi$ while $C_2=Q_2+\ell_2$, where $Q_2$ is a smooth conic, $\ell_2$ is a line tangent to $Q_2$ at $\xi$, and $\ell_2$ is also one of the two tangent lines of $C_1$ at $\xi$. We also require the configuration $C_1+C_2$ to be generic. Then, locally analytically $(C_1+C_2,\xi)$ is of type $Z_{15}$.

\item The curve $C_1$ has a node at $\xi$, and $C_2$ is the union of three generic lines through $\xi$ which are distinct from the two tangent lines of $C_1$. Then, locally analytically $(C_1+C_2,\xi)$ is of type $N_{16}$.

\end{enumerate}
The configuration $C_1,C_2$ is pictured in Figure~\ref{fig:Z_12-W_15-Z_1,0-N_16-singularity}. Let $\Sigma\in\{Z_{12},W_{15},Z_{15},N_{16}\}$ be one of the four cases above. If $L$ denotes a line in $\mathbb{P}^2$ not passing through $\xi$ and such that $C_1+C_2+L$ has only nodal singularities away from $\xi$, then we call the curve $C_1+C_2+L$ a \emph{$\Sigma$ arrangement}.
\end{lemma-definition}

\begin{proof}
The determination of the four singularities above is standard, but we write the proof for the reader's convenience. Consider the completed local ring $\widehat{\mathcal{O}}_{\mathbb{P}^2,\xi}\cong\mathbb{C}[[y,z]]$.

\medskip

\textbf{$Z_{12}$ arrangement.} Since $C_1$ has a node at $\xi$, we can choose analytic coordinates $y$ and $z$ such that locally analytically at $\xi$ the equation of $C_1$ is $yz=0$. Therefore, the cusp of $C_2$ is tangent at $\xi$ to, say, $z=0$. By assigning weights $\mathrm{wt}(y)=2,\mathrm{wt}(z)=3$, since the contact of $C_2$ with $\{z=0\}$ at $\xi$ is three, we have that
\[
C_2\colon z^2+\lambda y^3+(\text{terms of weighted degree \(>6\)})=0
\]
for some $\lambda\neq0$, which we can assume to be $1$ after rescaling $y$. In particular, locally analytically at $\xi$, the curve $C_1+C_2$ is described by the vanishing of
\[
f:=yz(y^3+z^2)+(\text{terms of weighted degree \(>11\)}).
\]
The weighted principal part of $C_1+C_2$ is therefore $f_0:=yz(y^3+z^2)$, whose weighted degree is $11$. By \cite[Ch.~I, Corollary~2.18]{GLS25} we have that the Milnor number $\mu(f)$ equals $\mu(f_0)$, which can be computed to be $12$. Therefore, we can apply Arnol'd's classification \cite{Arn76} (which is for singularities of modality at most $2$ or Milnor number at most $16$) to conclude that locally analytically at $\xi$ the curve $C_1+C_2$ has a $Z_{12}$ singularity.

\medskip

\textbf{$W_{15}$ arrangement.} We can choose projective coordinates so that the common tangent to the two cusps is $z=0$. The two cusp equations can be written as
\[
z^2+\lambda_i y^3+(\text{terms of weighted degree \(>6\)})=0,\qquad i=1,2,
\]
with weights $\mathrm{wt}(y)=2,\mathrm{wt}(z)=3$ and $\lambda_1,\lambda_2\neq0$. The condition $I_\xi(C_1,C_2)=6$ is equivalent to $\lambda_1\neq\lambda_2$. So, locally analytically at $\xi$, the curve $C_1+C_2$ is described by the vanishing of
\[
f:=z^4+(\lambda_1+\lambda_2)y^3z^2+\lambda_1\lambda_2y^6+(\text{terms of weighted degree \(>12\)}).
\]
After rescaling $y$, the weighted degree $12$ part becomes $y^6+a\,y^3z^2+z^4,\ a^2\neq4$. This has Milnor number $15$, and hence from \cite{Arn76} we can conclude that locally analytically at $\xi$ the curve $C_1+C_2$ has a $W_{15}$ singularity.

\medskip
 
\textbf{$Z_{15}$ arrangement.} The singularity at $\xi$ has four local branches. By choosing appropriate local analytic coordinates $(y,z)$ and after assigning weights $\mathrm{wt}(y)=1,\mathrm{wt}(z)=2$, we can write the germ $(C_1+C_2,\xi)$ as
\[
y\prod_{i=1}^3(z+\lambda_iy^2)+(\text{terms of weighted degree \(>7\)})=0,
\]
for generic constants $\lambda_1,\lambda_2,\lambda_3$. Expanding the product yields $\alpha y^6+\beta y^4z+\gamma y^2z^2+z^3$ for some generic coefficients $\alpha,\beta,\gamma$. After a linear change of coordinates in the variables $z$ and $y^2$ we can cancel the term $\beta y^4z$ and rescale $y^6$, so that the germ $(C_1+C_2,\xi)$ is described by
\[
y(y^6+ay^2z^2+z^3)+(\text{terms of weighted degree \(>7\)})=0,
\]
with generic $a$. This has Milnor number $15$, and hence locally analytically at $\xi$ the curve $C_1+C_2$ has a $Z_{15}$ singularity by the classification in \cite{Arn76}.

\medskip

\textbf{$N_{16}$ arrangement.} By choosing appropriate local analytic coordinates $(y,z)$, we can suppose that $C_1$ is locally described at $\xi$ by $yz=0$. Therefore, after assigning weights $\mathrm{wt}(y)=1,\mathrm{wt}(z)=1$, we can write the germ $(C_1+C_2,\xi)$ as
\[
yz\prod_{i=1}^3(a_iy+b_iz)+(\text{terms of weighted degree \(>5\)})=0,
\]
for generic constants $a_i,b_i$, $i=1,2,3$. This has Milnor number $16$, so we can use \cite{Arn76} to classify it. In conclusion, locally analytically at $\xi$ the curve $C_1+C_2$ has an $N_{16}$ singularity.
\end{proof}

The next lemma is a standard calculation, and it can be obtained by taking an appropriate log resolution of the pair $(\mathbb{P}^2,C_1+C_2+L)$.

\begin{lemma}
\label{lem:ldts-of-Sigma-arrangements}
Let $C_1+C_2+L$ be a $\Sigma$ arrangement. Then, the log canonical threshold of the pair $(\mathbb{P}^2,C_1+C_2+L)$ is given in the following table:
\begin{table}[h!]
\centering
\renewcommand{\arraystretch}{1.4}
\begin{tabular}{|c|c|c|c|c|}
\hline
$\Sigma$ & $Z_{12}$ & $W_{15}$ & $Z_{15}$ & $N_{16}$\\
\hline
$\mathrm{lct}$ & $\frac{5}{11}$ & $\frac{5}{12}$ & $\frac{3}{7}$ & $\frac{2}{5}$\\
\hline
\end{tabular}
\end{table}

\noindent In particular, $\left(\mathbb{P}^2,\frac{1}{2}(C_1+C_2+L)\right)$ and its $(\mathbb{Z}/2\mathbb{Z})^2$-cover $S_0\rightarrow\mathbb{P}^2$ branched along $C_1+C_2+L$ are not log canonical.
\end{lemma}

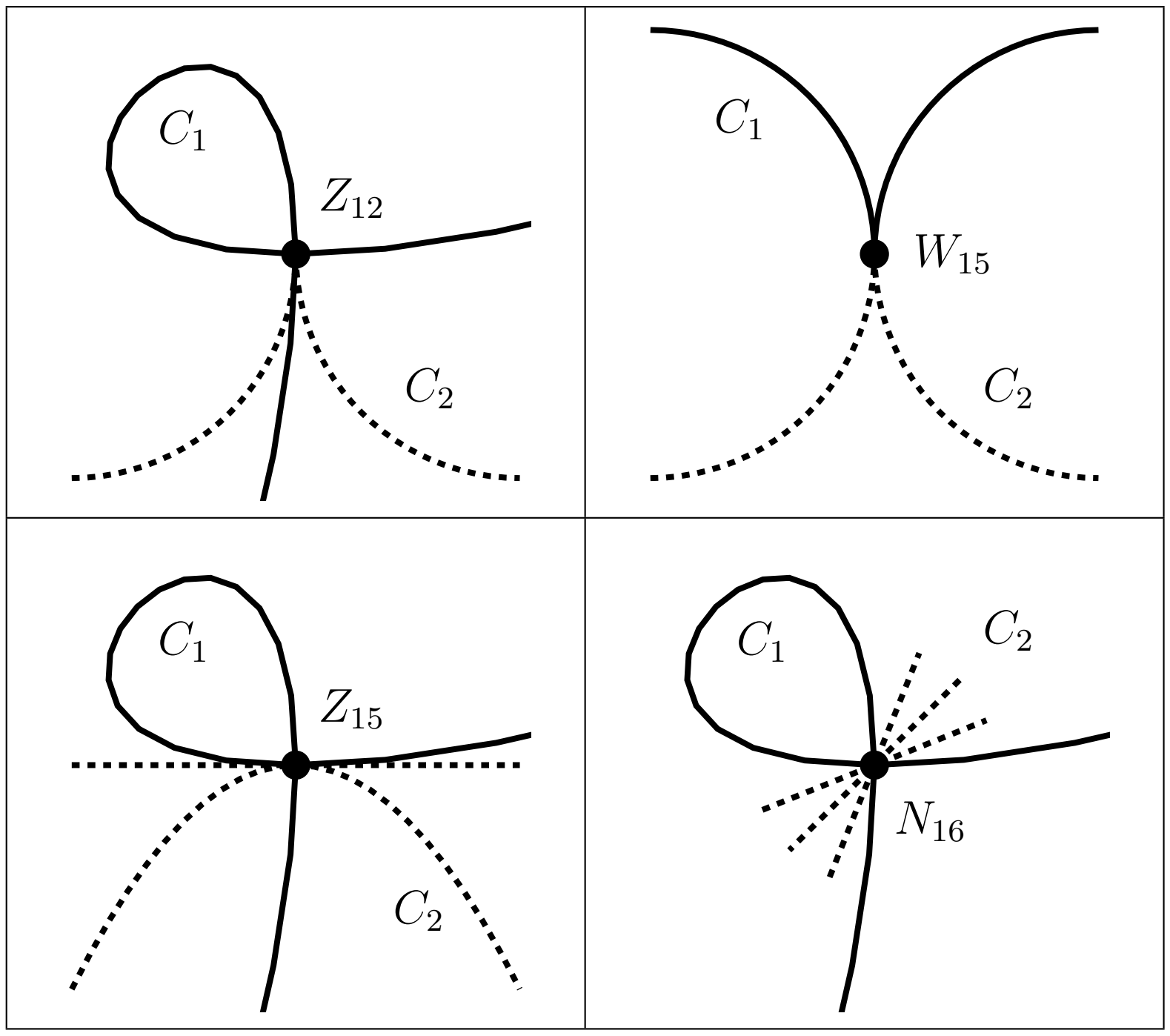
\begin{table}[hbtp]
\centering
\caption{The two plane cubics $C_1,C_2$ in Lemma--Definition~\ref{def:Z_12-W_15-Z_1,0-N_16-arrangements}.}
\label{fig:Z_12-W_15-Z_1,0-N_16-singularity}
\begin{tabular}{|>{\centering\arraybackslash}m{4.8cm}|>{\centering\arraybackslash}m{4.8cm}|}
\hline
\begin{tikzpicture}

\clip (-2.1,-2.2) rectangle (2.1,2.2);

\draw [shift={(-2.,0.)},line width=1.5pt,dash pattern=on 2pt off 2pt]  plot[domain=-1.5707963267948966:0.,variable=\t]({1.*2.*cos(\t r)+0.*2.*sin(\t r)},{0.*2.*cos(\t r)+1.*2.*sin(\t r)});
\draw [shift={(2.,0.)},line width=1.5pt,dash pattern=on 2pt off 2pt]  plot[domain=3.141592653589793:4.71238898038469,variable=\t]({1.*2.*cos(\t r)+0.*2.*sin(\t r)},{0.*2.*cos(\t r)+1.*2.*sin(\t r)});
\draw [line width=1.5pt]  plot[domain=-1.5:1.5,variable=\t]({(2^(1/2))*(\t*\t-1)+(2^(1/2))*(\t*\t*\t-\t)},{-(2^(1/2))*(\t*\t-1)+(2^(1/2))*(\t*\t*\t-\t)});
\draw (-1.0,1.1) node {$C_1$};
\draw (1.2,-1.2) node {$C_2$};
\draw (0.5,0.5) node {$Z_{12}$};

\draw [fill=black] (0,0) circle (3.5pt);

\end{tikzpicture}
&
\begin{tikzpicture}

\clip (-2.1,-2.2) rectangle (2.1,2.2);

\draw [shift={(-2.,0.)},line width=1.5pt,dash pattern=on 2pt off 2pt]  plot[domain=-1.5707963267948966:0.,variable=\t]({1.*2.*cos(\t r)+0.*2.*sin(\t r)},{0.*2.*cos(\t r)+1.*2.*sin(\t r)});
\draw [shift={(2.,0.)},line width=1.5pt,dash pattern=on 2pt off 2pt]  plot[domain=3.141592653589793:4.71238898038469,variable=\t]({1.*2.*cos(\t r)+0.*2.*sin(\t r)},{0.*2.*cos(\t r)+1.*2.*sin(\t r)});
\draw [shift={(-2.,0.)},line width=1.5pt]  plot[domain=0.:1.5707963267948966,variable=\t]({1.*2.*cos(\t r)+0.*2.*sin(\t r)},{0.*2.*cos(\t r)+1.*2.*sin(\t r)});
\draw [shift={(2.,0.)},line width=1.5pt]  plot[domain=1.5707963267948966:3.141592653589793,variable=\t]({1.*2.*cos(\t r)+0.*2.*sin(\t r)},{0.*2.*cos(\t r)+1.*2.*sin(\t r)});
\draw (-1.2,1.2) node {$C_1$};
\draw (1.2,-1.2) node {$C_2$};
\draw (0.7,0) node {$W_{15}$};

\draw [fill=black] (0.,0.) circle (3.5pt);

\end{tikzpicture}
\\
\hline
\begin{tikzpicture}

\clip (-2.1,-2.2) rectangle (2.1,2.2);

\draw [line width=1.5pt,dash pattern=on 2pt off 2pt] (-2.,0.)-- (2.,0.);
\draw [line width=1.5pt]  plot[domain=-1.5:1.5,variable=\t]({(2^(1/2))*(\t*\t-1)+(2^(1/2))*(\t*\t*\t-\t)},{-(2^(1/2))*(\t*\t-1)+(2^(1/2))*(\t*\t*\t-\t)});
\draw [line width=1.5pt,dash pattern=on 2pt off 2pt]  plot[domain=-2:2.,variable=\t]({\t},{-0.5*\t*\t});
\draw (1.1,-1.3) node {$C_2$};
\draw (-1.0,1.1) node {$C_1$};
\draw (0.5,0.5) node {$Z_{15}$};
\draw [fill=black] (0,0) circle (3.5pt);

\end{tikzpicture}
&
\begin{tikzpicture}

\clip (-2.1,-2.2) rectangle (2.1,2.2);

\draw [line width=1.5pt]  plot[domain=-1.5:1.5,variable=\t]({(2^(1/2))*(\t*\t-1)+(2^(1/2))*(\t*\t*\t-\t)},{-(2^(1/2))*(\t*\t-1)+(2^(1/2))*(\t*\t*\t-\t)});
\draw [line width=1.5pt,dash pattern=on 2pt off 2pt] (-0.4,-1)-- (0.4,1);
\draw [line width=1.5pt,dash pattern=on 2pt off 2pt] (-1,-0.4)-- (1,0.4);
\draw [line width=1.5pt,dash pattern=on 2pt off 2pt] (0.76,0.76)-- (-0.76,-0.76);
\draw (-1.0,1.1) node {$C_1$};
\draw (1.2,1.2) node {$C_2$};
\draw (0.5,-0.5) node {$N_{16}$};

\draw [fill=black] (0,0) circle (3.5pt);

\end{tikzpicture}
\\
\hline
\end{tabular}
\end{table}

\begin{remark}
From the previous lemma, we can see at once that the KSBA compactification $\overline{\mathrm{M}}_K$ differs from the GIT compactification of plane septic curves: by \cite{KL04}, if $\Sigma=Z_{12}$ then $C_1+C_2+L$ is GIT-stable because $\mathrm{lct}(\PP^2,C_1+C_2+L)=\frac{5}{11} > \frac{3}{7}$. If $\Sigma=Z_{15}$, then $C_1+C_2+L$ is strictly GIT-semistable because $\mathrm{lct}(\PP^2,C_1+C_2+L)= \frac{3}{7}$. On the other hand, for $\Sigma\in\{W_{15},N_{16}\}$, the curve $C_1+C_2+L$ is not GIT-semistable.
\end{remark}

\begin{proposition}
\label{prop:rationality-of-non-lc-limit}
Let $\Sigma\in\{Z_{12},W_{15},Z_{15},N_{16}\}$. Let $S_0$ be the $(\mathbb{Z}/2\mathbb{Z})^2$-cover of $\mathbb{P}^2$ branched along a $\Sigma$ arrangement $C_1+C_2+L$ as in Lemma--Definition~\ref{def:Z_12-W_15-Z_1,0-N_16-arrangements}. Then, $S_0$ is a rational surface.
\end{proposition}

\begin{proof}
Let $\xi\in C_1\cap C_2$ be the singular point of type $\Sigma$ and let $\ell\subseteq\mathbb{P}^2$ be a general line passing through $\xi$. Then, under the $(\mathbb{Z}/2\mathbb{Z})^2$-cover $\pi\colon S_0\rightarrow\mathbb{P}^2$, the curve $\widetilde{\ell}=\pi^*\ell$ is a rational curve. To prove this, let us write $\pi=\pi_1\circ\pi_2$, where $\pi_1\colon X\rightarrow\mathbb{P}^2$ is the double cover of $\mathbb{P}^2$ branched along $C_1+C_2$ and $\pi_2\colon S_0\rightarrow X$ is the double cover of $X$ with branch divisor $\pi^*L$. We have that $\pi_1^*\ell$ is a rational curve singular over $\xi$ that intersects $\pi_1^*L$ transversely in exactly two points. Hence, $\pi^*\ell=\pi_2^*\pi_1^*L$ is also a rational curve with a singularity over $\xi$. Therefore, the surface $S_0$ has a pencil of rational curves, which implies the claimed rationality.
\end{proof}

\section{Stable degenerations of Kunev surfaces}
\label{sec:stable-degen-of-Kunev-surfaces}

The objective of this section is to construct explicit one-parameter degenerations of Kunev surfaces whose central fiber is stable. These are constructed in the following way:
\begin{enumerate}

\item We construct one-parameter degenerations $\left(\mathcal{X},\frac{1}{2}\mathcal{B}\right)\rightarrow\Delta$ of the pairs given by $\left(\mathbb{P}^2,\frac{1}{2}(C_1+C_2+L)\right)$ such that the limit curve $B_0$ for $t=0$ has a $\Sigma$ singularity at the point $\xi=[1:0:0]$, $\Sigma\in\{Z_{12},W_{15},Z_{15},N_{16}\}$, and nodal singularities elsewhere. This is developed in Section~\ref{sec:non-lc-limits}.

\item We perform a weighted blow up of the central fiber of $\left(\mathcal{X},\frac{1}{2}\mathcal{B}\right)\rightarrow\Delta$ resolving the unique non-log canonical center, obtaining $\left(\mathcal{X}',\frac{1}{2}\mathcal{B}'\right)\rightarrow\Delta$. Then, we describe the geometry of the resulting exceptional surface $Y\subseteq\mathcal{X}'$ together with the curve $B_0'|_Y$. See Section~\ref{sec:blow-up-central-fiber}.

\item We describe the KSBA-stability of the new central fiber. In some cases, it is necessary to consider a birational modification $\mathcal{X}'\rightarrow\mathcal{X}''$ contracting certain curves in the central fiber to obtain the stable model $\left(\mathcal{X}'',\frac{1}{2}\mathcal{B}''\right)\rightarrow\Delta$. This is discussed in Section~\ref{sec:KSBA-stability-of-central-fiber-and-contraction}.

\end{enumerate}

Lastly, in Section~\ref{sec:description-of-the-cover-of-Y} we describe the geometry of the $(\mathbb{Z}/2\mathbb{Z})^2$-cover of the exceptional surface $Y$, which is the actual irreducible component of the stable degeneration of Kunev surfaces.

\subsection{One-parameter degenerations with non-log canonical central fiber}
\label{sec:non-lc-limits}

Consider a generic plane curve $C_1+C_2+L$, where $C_i=V(\gamma_i)$ is a cubic, $i=1,2$, and $L=V(\ell)$ is a line. We deform $C_1+C_2+L$ by appropriately introducing a parameter $t$ among the coefficients of $\gamma_1\cdot\gamma_2\cdot\ell$. The key definition is the following.

\begin{definition}
\label{def:Cstar-action}
Let $p,q,r$ be positive integers. For a monomial $m=x^iy^jz^k$, define the weight function
\[
\mathrm{wt}(x^iy^jz^k):=pj+qk-r.
\]
Given $t\in\mathbb{C}^*$, also let
\begin{displaymath}
t\star m:=\left\{ \begin{array}{ll}
t^{-\mathrm{wt}(m)}m&\textrm{if}~\mathrm{wt}(w)\in\mathbb{Z}_{\leq0}\\
m&\textrm{if}~\mathrm{wt}(w)\in\mathbb{Z}_{>0}.
\end{array} \right.
\end{displaymath}
By taking the linear extension of the above to $\mathbb{C}[x,y,z]$, we obtain an action
\[
\star\colon\mathbb{C}^*\times\mathbb{C}[x,y,z]\rightarrow\mathbb{C}[x,y,z].
\]
\end{definition}

\begin{definition}
\label{def:weight-functions}
Given $\Sigma\in\{Z_{12},W_{15},Z_{15},N_{16}\}$, we define two weight functions $\mathrm{wt}_1$ and $\mathrm{wt}_2$ as follows:

\begin{itemize}

\item If $\Sigma=Z_{12}$, let $\mathrm{wt}_1(x^iy^jz^k)=2j+3k-5$, $\mathrm{wt}_2(x^iy^jz^k)=2j+3k-6$;

\item If $\Sigma=W_{15}$, let $\mathrm{wt}_1(x^iy^jz^k)=2j+3k-6=\mathrm{wt}_2(x^iy^jz^k)$;

\item If $\Sigma=Z_{15}$, let $\mathrm{wt}_1(x^iy^jz^k)=j+2k-3$, $\mathrm{wt}_2(x^iy^jz^k)=j+2k-4$;

\item If $\Sigma=N_{16}$, let $\mathrm{wt}_1(x^iy^jz^k)=j+k-2$, $\mathrm{wt}_2(x^iy^jz^k)=j+k-3$.

\end{itemize}
We denote by $\star_1$ and $\star_2$ the corresponding $\mathbb{C}^*$-actions on $\mathbb{C}[x,y,z]$ in the sense of Definition~\ref{def:Cstar-action}.
\end{definition}

\begin{lemma}
\label{lem:Sigma-singularity-at-the-limit}
Let $\Sigma\in\{Z_{12},W_{15},Z_{15},N_{16}\}$ and consider the corresponding weight functions $\mathrm{wt}_1,\mathrm{wt}_2$ as in Definition~\ref{def:weight-functions}. Let $F_1,F_2\in\mathbb{C}[x,y,z]$ be two degree three polynomials with generic coefficients. Then, the plane curve defined by
\[
\left(\lim_{t\to0}t\star_1F_1\right)\cdot\left(\lim_{t\to0}t\star_2F_2\right)=0
\]
has a $\Sigma$ singularity at $\xi=[1:0:0]$ and it is at worst nodal everywhere else.
\end{lemma}

\begin{proof}
For a degree $3$ monomial in $\mathbb{C}[x,y,z]$, the corresponding values of $\mathrm{wt}_1$ and $\mathrm{wt}_2$ are listed in Table~\ref{tbl:weights-of-monomials-degree-3-in-P2-depending-on-Sigma}. With these values, we can compute $t\star_1F_1$, $t\star_2F_2$, and the corresponding limits for $t\to0$ for each $\Sigma$ individually.

\begin{itemize}

\item For $\Sigma=Z_{12}$, we have that
\begin{align*}
t\star_1F_1&=t^5x^3+t^3x^2y+t^2x^2z+txy^2+xyz+xz^2+y^3+y^2z+yz^2+z^3,\\
t\star_2F_2&=t^6x^3+t^4x^2y+t^3x^2z+t^2xy^2+txyz+xz^2+y^3+y^2z+yz^2+z^3,
\end{align*}
where the coefficients are intended to be generic. Hence,
\begin{align*}
\lim_{t\to0}t\star_1F_1&=xyz+xz^2+y^3+y^2z+yz^2+z^3,\\
\lim_{t\to0}t\star_2F_2&=xz^2+y^3+y^2z+yz^2+z^3.
\end{align*}
The corresponding plane curves are singular at $[1:0:0]$, where they have a node and a cusp, respectively. Moreover, one of the two tangents to the node coincides with the double tangent line of the cusp. Hence, the plane curve
\begin{equation}
\label{eq:limit-curve-Z12-case}
\left(\lim_{t\to0}t\star_1F_1\right)\cdot\left(\lim_{t\to0}t\star_2F_2\right)=0
\end{equation}
has a $Z_{12}$ singularity at $[1:0:0]$ by Lemma--Definition~\ref{def:Z_12-W_15-Z_1,0-N_16-arrangements}~(1).

Let $\ell_i$ be the linear system generated by the monomials of $\lim_{t\to0}t\star_iF_i$, $i=1,2$. We have that $\ell_i$ has $\xi$ as unique base point. Hence, by Bertini's theorem \cite[Ch.~1]{GH94}, we have that the generic curve $C_1\in\ell_1$ is smooth away from $\xi$ and the restricted linear system $\ell_2|_{C_1}$ is base point free away from $\xi$, and hence it consists of reduced points. As the coefficients are generic, this implies that away from the $Z_{12}$ singularity the curve in \eqref{eq:limit-curve-Z12-case} has at worst nodes.

\item For $\Sigma=W_{15}$, we have that
\[
t\star_1F_1=t^6x^3+t^4x^2y+t^3x^2z+t^2xy^2+txyz+xz^2+y^3+y^2z+yz^2+z^3.
\]
By taking the limit $t\to0$, we obtain that
\[
\lim_{t\to0}t\star F_1=xz^2+y^3+y^2z+yz^2+z^3.
\]
The corresponding plane curves are singular at $[1:0:0]$, where it has a cusp with tangent $z=0$. An analogous calculation holds for $t\star_2F_2$, from which we obtain that also $\lim_{t\to0}t\star F_2$ has a cusp at $[1:0:0]$ with tangent $z=0$. By the genericity of the coefficients of $F_1,F_2$, we can conclude that the plane curve
\[
\left(\lim_{t\to0}t\star F_1\right)\cdot\left(\lim_{t\to0}t\star F_2\right)=0
\]
has a $W_{15}$ singularity at $[1:0:0]$ by Lemma--Definition~\ref{def:Z_12-W_15-Z_1,0-N_16-arrangements}~(2), with nodal singularities elsewhere (this last assertion can be justified as in the $Z_{12}$ case).

\item For $\Sigma=Z_{15}$, we have that
\begin{align*}
t\star_1F_1&=t^3x^3+t^2x^2y+tx^2z+txy^2+xyz+xz^2+y^3+y^2z+yz^2+z^3,\\
t\star_2F_2&=t^4x^3+t^3x^2y+t^2x^2z+t^2xy^2+txyz+xz^2+ty^3+y^2z+yz^2+z^3.
\end{align*}
By taking the limit $t\to0$, we obtain that
\begin{align*}
\lim_{t\to0}t\star_1F_1&=xyz+xz^2+y^3+y^2z+yz^2+z^3\\
&=z^2(x+y+z)+yz(x+y)+y^3,\\
\lim_{t\to0}t\star_2F_2&=xz^2+y^2z+yz^2+z^3\\
&=z^2(x+y+z)+y^2z.
\end{align*}
The first curve is a cubic nodal at $\xi$ with tangent directions $z=0$ and $y+z=0$. The second curve is the union of a smooth conic together with the line $z=0$ such that the line is tangent to the conic at $\xi$. Hence, we can conclude that the plane curve
\[
\left(\lim_{t\to0}t\star_1F_1\right)\cdot\left(\lim_{t\to0}t\star_2F_2\right)=0
\]
has a $Z_{15}$ singularity at $\xi$ by Lemma--Definition~\ref{def:Z_12-W_15-Z_1,0-N_16-arrangements}~(3). To prove that the only other singular points are nodes located at the intersection points between the two curves, also here we use Bertini's theorem as in the $Z_{12}$ case, but with the following adaptation. In this $Z_{15}$ case, the linear system $\ell_2$ has the line $z=0$ as its base locus. At the same time, this line intersects $C_1$ only at the point $\xi$. So, we can still conclude that away from $\xi$ the singularities are nodes.

\item For $\Sigma=N_{16}$, we have that
\begin{align*}
t\star_1F_1&=t^2x^3+tx^2y+tx^2z+xy^2+xyz+xz^2+y^3+y^2z+yz^2+z^3,\\
t\star_2F_2&=t^3x^3+t^2x^2y+t^2x^2z+txy^2+txyz+txz^2+y^3+y^2z+yz^2+z^3.
\end{align*}
By taking the limit $t\to0$, we obtain that
\begin{align*}
\lim_{t\to0}t\star_1F_1&=xy^2+xyz+xz^2+y^3+y^2z+yz^2+z^3,\\
\lim_{t\to0}t\star_2F_2&=y^3+y^2z+yz^2+z^3.
\end{align*}
The corresponding plane curves are singular at $[1:0:0]$ giving rise to an $N_{16}$ singularity. By the genericity of the coefficients of $F_1,F_2$, we can conclude that the plane curve
\[
\left(\lim_{t\to0}t\star F_1\right)\cdot\left(\lim_{t\to0}t\star F_2\right)=0
\]
has an $N_{16}$ singularity at $[1:0:0]$, with nodal singularities elsewhere (this last assertion can be justified as in the $Z_{12}$ case).\qedhere

\end{itemize}

\end{proof}

\begin{table}
\centering
\caption{The monomials $x^iy^jz^k$ of degree $3$ and the corresponding values of the weight functions $\mathrm{wt}_1$ and $\mathrm{wt}_2$ depending on the singularity $\Sigma\in\{Z_{12},W_{15},Z_{15},N_{16}\}$.}
\label{tbl:weights-of-monomials-degree-3-in-P2-depending-on-Sigma}
\renewcommand{\arraystretch}{1.4}
\begin{tabular}{|c|c|c|c|c|c|c|c|c|c|c|c|}
\hline
\multicolumn{2}{|c|}{} &  $x^3$ & $x^2y$ & $x^2z$ & $xy^2$ & $xyz$ & $xz^2$ & $y^3$ & $y^2z$ & $yz^2$ & $z^3$ \\
\hline
\multirow{2}{*}{$Z_{12}$}& $\mathrm{wt}_1$ & $-5$ & $-3$ & $-2$ & $-1$ & $0$ & $1$ & $1$ & $2$ & $3$ & $4$ \\ \cline{2-12}
& $\mathrm{wt}_2$ & $-6$ & $-4$ & $-3$ & $-2$ & $-1$ & $0$ & $0$ & $1$ & $2$ & $3$ \\
\hline
$W_{15}$ & $\mathrm{wt}_{1,2}$ & $-6$ & $-4$ & $-3$ & $-2$ & $-1$ & $0$ & $0$ & $1$ & $2$ & $3$ \\
\hline
\multirow{2}{*}{$Z_{15}$}& $\mathrm{wt}_1$ & $-3$ & $-2$ & $-1$ & $-1$ & $0$ & $1$ & $0$ & $1$ & $2$ & $3$ \\ \cline{2-12}
& $\mathrm{wt}_2$ & $-4$ & $-3$ & $-2$ & $-2$ & $-1$ & $0$ & $-1$ & $0$ & $1$ & $2$ \\
\hline
\multirow{2}{*}{$N_{16}$}& $\mathrm{wt}_1$ & $-2$ & $-1$ & $-1$ & $0$ & $0$ & $0$ & $1$ & $1$ & $1$ & $1$ \\ \cline{2-12}
& $\mathrm{wt}_2$ & $-3$ & $-2$ & $-2$ & $-1$ & $-1$ & $-1$ & $0$ & $0$ & $0$ & $0$ \\
\hline
\end{tabular}
\end{table}

\begin{definition}
\label{def:1-parameter-degeneration-Sigma}
Let $\Delta=\mathrm{Spec}(R)$ with $R=\mathbb{C}[[t]]$ and $\mathcal{X}=\Delta\times\mathbb{P}^2$. Let us consider general homogeneous polynomials $F_1,F_2,\ell\in\mathbb{C}[x,y,z]$ with $\deg(F_1)=\deg(F_2)=3$, $\deg(\ell)=1$, and define $\mathcal{B}\subseteq\mathcal{X}$ to be the divisor
\[
\mathcal{B}\colon(t\star_1F_1)\cdot(t\star_2F_2)\cdot\ell=0.
\]
Then, $\left(\mathcal{X},\frac{1}{2}\mathcal{B}\right)\rightarrow\Delta$ is a family of pairs such that, for generic $t$, the fiber $\left(X_t=\mathbb{P}^2,\frac{1}{2}B_t\right)$ has $B_t$ given by the union of two general cubics and a general line, which gives a stable pair. On the other hand, for $t=0$, $\left(\mathbb{P}^2,\frac{1}{2}B_0\right)$ is not log canonical. This is because by Lemma~\ref{lem:Sigma-singularity-at-the-limit} the curve $B_0$ is a $\Sigma$ arrangement, and its log canonical threshold, given in Lemma~\ref{lem:ldts-of-Sigma-arrangements}, is strictly less than $\frac{1}{2}$ for any choice of $\Sigma$. Hence, the limit pair $\left(\mathbb{P}^2,\frac{1}{2}B_0\right)$ is not log canonical.

By taking an appropriate $(\mathbb{Z}/2\mathbb{Z})^2$-cover $\mathcal{S}\rightarrow\mathcal{X}$ branched along $\mathcal{B}$ we obtain the corresponding one parameter degeneration of Kunev surfaces $\mathcal{S}\rightarrow\Delta$ with non-log canonical central fiber $S_0$.
\end{definition}

\subsection{Blow up of the central fiber}
\label{sec:blow-up-central-fiber}

\begin{definition}
\label{def:weighted-blow-up-of-family-Sigma}
Let $\left(\mathcal{X},\frac{1}{2}\mathcal{B}\right)\rightarrow\Delta$ as in Definition~\ref{def:1-parameter-degeneration-Sigma}. Then, the curve $B_0$ has a $\Sigma$ singularity at $(t,[x:y:z])=(0,[1:0:0])$. By considering the affine patch $x\neq0$, we have affine coordinates $(t,y,z)$. Let $\mathcal{X}'\rightarrow\mathcal{X}$ be the weighted blow up with weights $(1,p,q)$ at $(t,y,z)=(0,0,0)$, where $p,q$ are as prescribed in the proof of Lemma~\ref{lem:ldts-of-Sigma-arrangements}. Denote by $\mathcal{B}'\subseteq\mathcal{X}'$ be the strict transform of $\mathcal{B}$, $Y\subseteq\mathcal{X}'$ the exceptional divisor, and $Z\subseteq\mathcal{X}'$ the strict transform of $X_0\subseteq\mathcal{X}$. The new central fiber $\left(X_0',\frac{1}{2}B_0'\right)$ of $(\mathcal{X}',\frac{1}{2}\mathcal{B}')\rightarrow\Delta$ consists of two irreducible components:
\begin{enumerate}

\item $\left(Y,T+\frac{1}{2}B_0'|_Y\right)$, where $Y\cong\mathbb{P}(1,p,q)$ with coordinates $[\tau:\alpha:\beta]$, $T=V(\tau)$, and $B_0'|_Y$ is a degree $d_\Sigma$ curve (see again the proof of Lemma~\ref{lem:ldts-of-Sigma-arrangements} for the definition of $d_\Sigma$) given by

\begin{itemize}

\item $\Sigma=Z_{12}$:
\begin{align*}
(&\tau^5+\tau^3\alpha+\tau^2\beta+\tau\alpha^2+\alpha\beta)\\
\cdot(&\tau^6+\tau^4\alpha+\tau^3\beta+\tau^2\alpha^2+\tau \alpha\beta+\alpha^3+\beta^2)=0,
\end{align*}
To obtain such equation, from the product $(t\star_1F_1)\cdot(t\star_2F_2)\cdot(x+y+z)$, we dehomogeneize with respect to $x$, assign weights $(1,p,q)=(1,2,3)$ to $(t,y,z)$, take the lowest weight part, and substitute $\tau,\alpha,\beta$ to the variables $t,y,z$, respectively (the same procedure is applied in the other cases). The curve $B_0'|_Y$ intersects $T$ at the points $[0:1:0],[0:0:1],[0:1:-1]$.

\item $\Sigma=W_{15}$:
\begin{align*}
(&\tau^6+\tau^4\alpha+\tau^3\beta+\tau^2\alpha^2+\tau\alpha\beta+\beta^2+\alpha^3)\\
\cdot(&\tau^6+\tau^4\alpha+\tau^3\beta+\tau^2\alpha^2+\tau\alpha\beta+\beta^2+\alpha^3)=0,
\end{align*}
where the coefficients are general. By setting $\tau=0$, we obtain
\[
(b_1\beta^2+a_1\alpha^3)\cdot(b_2\beta^2+a_2\alpha^3)=0,
\]
(we now include explicit general coefficients) whose solutions give the two distinct points $[0:b_1:-a_1]$ and $[0:b_2:-a_2]$.

\item $\Sigma=Z_{15}$:
\begin{align*}
&(\tau^3+\tau^2\alpha+\tau\beta+\tau\alpha^2+\alpha\beta+\alpha^3)\\
\cdot&(\tau^4+\tau^3\alpha+\tau^2\alpha^2+\tau^2\beta+\tau\alpha^3+\tau\alpha\beta+\alpha^2\beta+\beta^2)=0.
\end{align*}
We note that in the second factor the monomial $\alpha^4$ is missing. However, as we will see the proof of Lemma~\ref{lem:moduli-for-curves-in-P1pq}, given a general quartic curve in $\mathbb{P}(1,1,2)$, we can assume that its defining polynomial does not have the monomial $\alpha^4$ after applying an appropriate automorphism. The curve $B_0'|_Y$ intersects $T$ at the points $[0:1:0],[0:0:1],[0:1:c],[0:1:c']$.

\item $\Sigma=N_{16}$:
\begin{align*}
(&\tau^2+\tau\alpha+\tau\beta+\alpha^2+\alpha\beta+\beta^2)\\
\cdot(&\tau^3+\tau^2\alpha+\tau^2\beta+\tau\alpha^2+\tau\alpha\beta+\tau\beta^2+\alpha^3+\alpha^2\beta+\alpha \beta^2+\beta^3)=0.
\end{align*}
The curve $B_0'|_Y$ intersects $T$ at five general points.

\end{itemize}

\item $\left(Z,E+\frac{1}{2}B_0'|_Z\right)$, where $Z$ is the blow up of $\mathbb{P}^2$ at $[1:0:0]$ with weights $(p,q)$, $E\subseteq Z$ is the exceptional divisor, and $B_0'|_Z$ is the strict transform of the singular curve $B_0\subseteq\mathbb{P}^2$.

\end{enumerate}
We point out that in each case the curve $B_0'|_Y$ is reducible, and it is given by the union of two generic curves of degrees $(d_\Sigma^{(1)},d_\Sigma^{(2)})$ given by
\begin{table}[h!]
\centering
\renewcommand{\arraystretch}{1.4}
\begin{tabular}{|c|c|c|c|c|}
\hline
$\Sigma$ & $Z_{12}$ & $W_{15}$ & $Z_{15}$ & $N_{16}$\\
\hline
$(d_\Sigma^{(1)},d_\Sigma^{(2)})$ & $(5,6)$ & $(6,6)$ & $(3,4)$ & $(2,3)$\\
\hline
\end{tabular}
\end{table}

In the cases $Z_{12},W_{15},N_{16}$ the two curves of degree $d_\Sigma^{(1)}$ and $d_\Sigma^{(2)}$ recover the generic degree $d_\Sigma^{(1)}$ and $d_\Sigma^{(2)}$ curves in $\mathbb{P}(1,p,q)$. This does not hold for $Z_{15}$: the degree $d_{Z_{15}}^{(2)}=4$ component misses the monomial $\alpha^4$.

The surfaces $Y$ and $Z$ are glued transversely along $Y\supseteq T\cong\mathbb{P}^1\cong E\subseteq Z$. We summarize this information in Figure~\ref{fig:slc-replacement-Sigma}.
\end{definition}

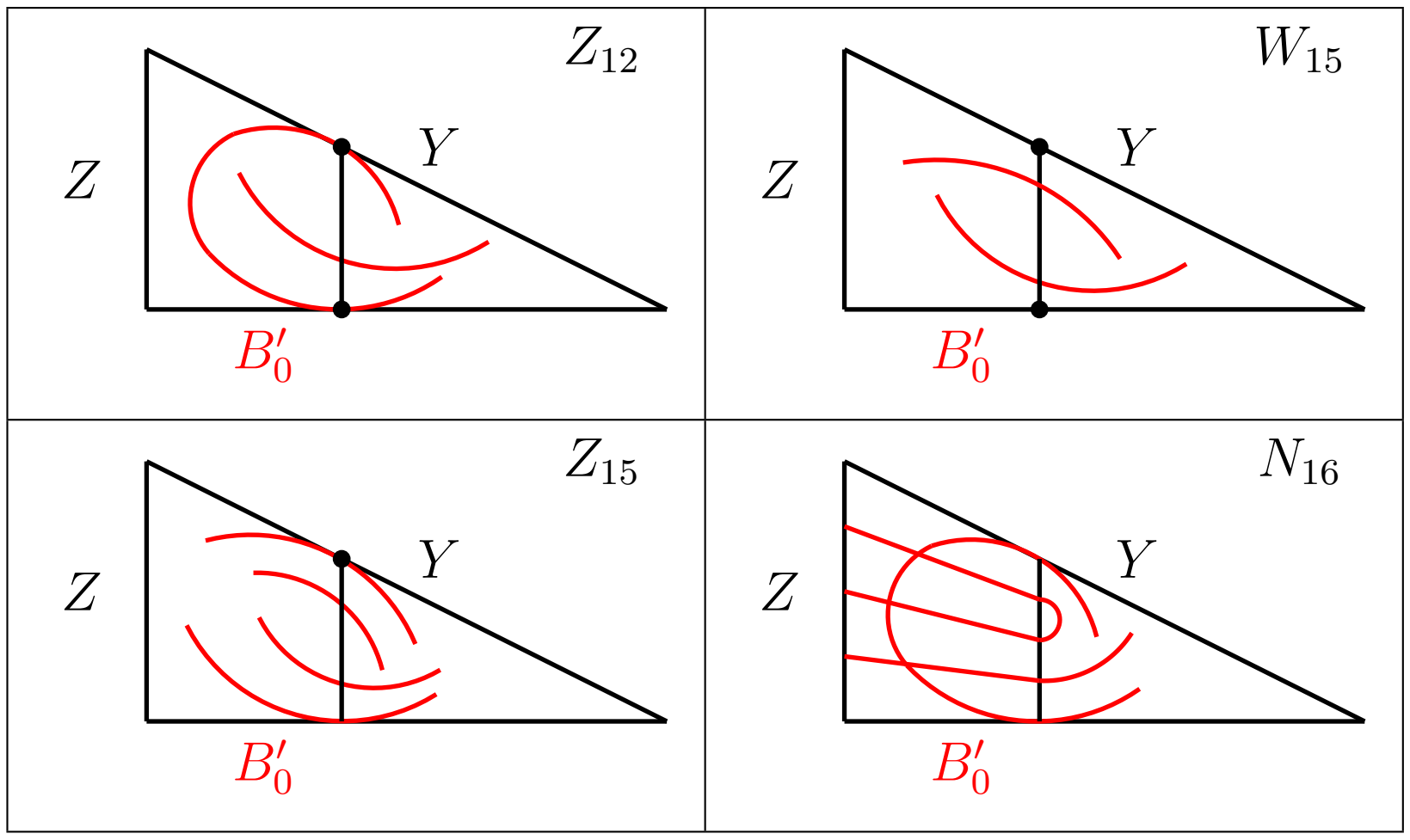
\begin{table}[hbtp]
\centering
\caption{The pair $\left(X_0',\frac{1}{2}B_0'\right)$ in Definition~\ref{def:weighted-blow-up-of-family-Sigma}.}
\label{fig:slc-replacement-Sigma}
\begin{tabular}{|>{\centering\arraybackslash}m{5.0cm}|>{\centering\arraybackslash}m{5.0cm}|}
\hline
\begin{tikzpicture}[scale=0.5]
\draw [line width=1.pt] (0.,4.)-- (0.,0.);
\draw [line width=1.pt] (0.,0.)-- (8.,0.);
\draw [line width=1.pt] (0.,4.)-- (8.,0.);
\draw [shift={(1.954405286343612,0.8105638766519823)},line width=1.pt,color=red]  plot[domain=0.2489032652766453:1.8851902781320244,variable=\t]({1.*1.9868222668621776*cos(\t r)+0.*1.9868222668621776*sin(\t r)},{0.*1.9868222668621776*cos(\t r)+1.*1.9868222668621776*sin(\t r)});
\draw [shift={(2.9581931411022007,2.7503651254052213)},line width=1.pt,color=red]  plot[domain=3.89926835953523:5.325066101674361,variable=\t]({1.*2.7506828491296806*cos(\t r)+0.*2.7506828491296806*sin(\t r)},{0.*2.7506828491296806*cos(\t r)+1.*2.7506828491296806*sin(\t r)});
\draw [shift={(1.8476383763837634,1.6359225092250922)},line width=1.pt,color=red]  plot[domain=2.0159314212059316:3.859936834018558,variable=\t]({1.*1.1789646421972833*cos(\t r)+0.*1.1789646421972833*sin(\t r)},{0.*1.1789646421972833*cos(\t r)+1.*1.1789646421972833*sin(\t r)});
\draw [shift={(3.823853329878207,3.3228271572946357)},line width=1.pt,color=red]  plot[domain=3.612171707888009:5.273937488941906,variable=\t]({1.*2.697001499477489*cos(\t r)+0.*2.697001499477489*sin(\t r)},{0.*2.697001499477489*cos(\t r)+1.*2.697001499477489*sin(\t r)});

\draw [fill=black] (3.,0.) circle (3.5pt);
\draw [fill=black] (3.,2.5) circle (3.5pt);

\draw [line width=1.pt] (3.,2.5)-- (3.,0.);
\draw (4.5,2.5) node {$Y$};
\draw (-1,2) node {$Z$};
\draw (1.8,-0.7) node {${\color{red}B_0'}$};
\draw (7.0,4.0) node {$Z_{12}$};
\end{tikzpicture}
&
\begin{tikzpicture}[scale=0.5]

\draw [line width=1.pt] (0.,4.)-- (0.,0.);
\draw [line width=1.pt] (0.,0.)-- (8.,0.);
\draw [line width=1.pt] (0.,4.)-- (8.,0.);
\draw [shift={(3.8238533298782067,2.982827157294636)},line width=1.pt,color=red]  plot[domain=3.6121717078880096:5.273937488941906,variable=\t]({1.*2.697001499477489*cos(\t r)+0.*2.697001499477489*sin(\t r)},{0.*2.697001499477489*cos(\t r)+1.*2.697001499477489*sin(\t r)});
\draw [line width=1.pt] (3.,2.5)-- (3.,0.);
\draw [shift={(1.423267326732674,-1.0678967085897775)},line width=1.pt,color=red]  plot[domain=0.5806113925115347:1.7267560584801565,variable=\t]({1.*3.3687839346370425*cos(\t r)+0.*3.3687839346370425*sin(\t r)},{0.*3.3687839346370425*cos(\t r)+1.*3.3687839346370425*sin(\t r)});

\draw [fill=black] (3.,0.) circle (3.5pt);
\draw [fill=black] (3.,2.5) circle (3.5pt);

\draw (4.5,2.5) node {$Y$};
\draw (-1,2) node {$Z$};
\draw (1.8,-0.7) node {${\color{red}B_0'}$};
\draw (7.0,4.0) node {$W_{15}$};
\end{tikzpicture}
\\
\hline
\begin{tikzpicture}[scale=0.5]
\draw [line width=1.pt] (0.,4.)-- (0.,0.);
\draw [line width=1.pt] (0.,0.)-- (8.,0.);
\draw [line width=1.pt] (0.,4.)-- (8.,0.);
\draw [shift={(3.5,2.5)},line width=1.pt,color=red]  
  plot[domain=0.2489032652766453:1.8851902781320244,variable=\t]
  (
    {1.9868222668621776*cos(\t r)*cos(-45) - (-1.9868222668621776*sin(\t r))*sin(-45)},
    {1.9868222668621776*cos(\t r)*sin(-45) + (-1.9868222668621776*sin(\t r))*cos(-45)}
  );
\draw [shift={(1.7,0.3)},line width=1.pt,color=red]  plot[domain=0.2489032652766453:1.6,variable=\t]({1.*1.9868222668621776*cos(\t r)+0.*1.9868222668621776*sin(\t r)},{0.*1.9868222668621776*cos(\t r)+1.*1.9868222668621776*sin(\t r)});
\draw [shift={(1.6,0.12)},line width=1.pt,color=red]  
  plot[domain=3.89926835953523:5.325066101674361,variable=\t]
  (
    {2.7506828491296806*cos(\t r)*cos(-32) - (-2.7506828491296806*sin(\t r))*sin(-32)},
    {2.7506828491296806*cos(\t r)*sin(-32) + (-2.7506828491296806*sin(\t r))*cos(-32)}
  );
\draw [shift={(3.02,2.7)},line width=1.pt,color=red]  plot[domain=3.612171707888009:5.273937488941906,variable=\t]({1.*2.697001499477489*cos(\t r)+0.*2.697001499477489*sin(\t r)},{0.*2.697001499477489*cos(\t r)+1.*2.697001499477489*sin(\t r)});

\draw [fill=black] (3.,2.5) circle (3.5pt);

\draw [line width=1.pt] (3.,2.5)-- (3.,0.);
\draw (4.5,2.5) node {$Y$};
\draw (-1,2) node {$Z$};
\draw (1.8,-0.7) node {${\color{red}B_0'}$};
\draw (7.0,4.0) node {$Z_{15}$};
\end{tikzpicture}
&
\begin{tikzpicture}[scale=0.5]
\draw [line width=1.pt] (0.,4.)-- (0.,0.);
\draw [line width=1.pt] (0.,0.)-- (8.,0.);
\draw [line width=1.pt] (0.,4.)-- (8.,0.);
\draw [shift={(1.954405286343612,0.8105638766519823)},line width=1.pt,color=red]  plot[domain=0.2489032652766453:1.8851902781320244,variable=\t]({1.*1.9868222668621776*cos(\t r)+0.*1.9868222668621776*sin(\t r)},{0.*1.9868222668621776*cos(\t r)+1.*1.9868222668621776*sin(\t r)});
\draw [shift={(2.9581931411022007,2.7503651254052213)},line width=1.pt,color=red]  plot[domain=3.89926835953523:5.325066101674361,variable=\t]({1.*2.7506828491296806*cos(\t r)+0.*2.7506828491296806*sin(\t r)},{0.*2.7506828491296806*cos(\t r)+1.*2.7506828491296806*sin(\t r)});
\draw [shift={(1.8476383763837634,1.6359225092250922)},line width=1.pt,color=red]  plot[domain=2.0159314212059316:3.859936834018558,variable=\t]({1.*1.1789646421972833*cos(\t r)+0.*1.1789646421972833*sin(\t r)},{0.*1.1789646421972833*cos(\t r)+1.*1.1789646421972833*sin(\t r)});
\draw [line width=1.pt] (3.,2.5)-- (3.,0.);
\draw [line width=1.pt,color=red] (0.,3.)-- (3.,1.875);
\draw [line width=1.pt,color=red] (0.,2.)-- (3.,1.25);
\draw [line width=1.pt,color=red] (0.,1.)-- (3.,0.625);
\draw [shift={(3.,1.5625)},line width=1.pt,color=red]  plot[domain=-1.5707963267948966:1.5707963267948966,variable=\t]({1.*0.3125*cos(\t r)+0.*0.3125*sin(\t r)},{0.*0.3125*cos(\t r)+1.*0.3125*sin(\t r)});
\draw [shift={(3.070517621145374,2.227962555066081)},line width=1.pt,color=red]  plot[domain=4.668425269112197:5.711618113988559,variable=\t]({1.*1.6045129129539533*cos(\t r)+0.*1.6045129129539533*sin(\t r)},{0.*1.6045129129539533*cos(\t r)+1.*1.6045129129539533*sin(\t r)});

\draw (4.5,2.5) node {$Y$};
\draw (-1,2) node {$Z$};
\draw (1.8,-0.7) node {${\color{red}B_0'}$};
\draw (7.0,4.0) node {$N_{16}$};
\end{tikzpicture}
\\
\hline
\end{tabular}
\end{table}

Next, we study the singularities and count the moduli for the pairs $\left(Y,T+\frac{1}{2}B_0'|_Y\right)$ for $Z_{12},W_{15},Z_{15}$. We omit $N_{16}$ from our discussion because $N_{16}$ will be dealt with in a different way (see Section~\ref{sec:stable-non-normal-K3-limits-from-N16}).

\begin{lemma}
\label{lem:moduli-for-curves-in-P1pq}
Let $\Sigma$ be one of the singularities $\{Z_{12},W_{15},Z_{15}\}$ and let us consider the corresponding pair $\left(Y,T+\frac{1}{2}B_0'|_Y\right)$ as in Definition~\ref{def:weighted-blow-up-of-family-Sigma}, where $Y\cong\mathbb{P}(1,p,q)$ and $B_0'|_Y$ is a reducible curve of degree $d_\Sigma$ with irreducible components of degrees $(d_\Sigma^{(1)},d_\Sigma^{(2)})$. Then,
\begin{enumerate}

\item The generic curve $B_0'|_Y$ has at worst nodal singularities away from the curve $\tau=0$.

\item The isomorphism class of the pairs $\left(Y,T+\frac{1}{2}B_0'|_Y\right)$ depends on the following parameters:
\begin{table}[h!]
\centering
\renewcommand{\arraystretch}{1.4}
\begin{tabular}{|c|c|c|c|}
\hline
$\Sigma$ & $Z_{12}$ & $W_{15}$ & $Z_{15}$ \\
\hline
$\mathrm{Moduli~for} \left(Y,T+\frac{1}{2}B_0'|_Y\right)$ & $5$ & $7$ & $7$ \\
\hline
\end{tabular}
\end{table}

\end{enumerate}
\end{lemma}

\begin{proof}
By the explicit description in Definition~\ref{def:weighted-blow-up-of-family-Sigma}, we can write $B_0'|_Y=V(f_\Sigma^{(1)}f_\Sigma^{(2)})$, where $f_\Sigma^{(1)},f_\Sigma^{(2)}$ are the polynomials defining the two components of $B_0'|_Y$ of degrees $d_\Sigma^{(1)},d_\Sigma^{(2)}$, respectively. After restricting to the smooth affine patch $\tau\neq0$, let $\ell_\Sigma^{(i)}$ be the linear system generated by the monomials of $f_\Sigma^{(i)}(1,\alpha,\beta)$, $i=1,2$. We have that $\ell_\Sigma^{(1)}$ and $\ell_\Sigma^{(2)}$ are base point free, so by Bertini's theorem we have that the generic curve $C_1\in\ell_\Sigma^{(1)}$ is smooth, and the restricted linear system $\ell_\Sigma^{(2)}|_{C_1}$ is base point free. In conclusion, for generic coefficients, $B_0'|_Z$ has at worst nodal singularities in the affine patch $\tau\neq0$. This proves (1).

To prove (2), let us first consider $\Sigma\in\{Z_{12},W_{15}\}$. Let
\[
\mathbb{P}(\mathbb{C}[\tau,\alpha,\beta]_{d_\Sigma^{(1)}}),\qquad\mathbb{P}(\mathbb{C}[\tau,\alpha,\beta]_{d_\Sigma^{(2)}}),\qquad\mathbb{P}(\mathbb{C}[\tau,\alpha,\beta]_{d_\Sigma})
\]
be the parameter spaces of curves in $\mathbb{P}(1,p,q)$ of degree $d_\Sigma^{(1)},d_\Sigma^{(2)},d_\Sigma$, respectively. Now, let $U\subseteq\mathbb{P}(\mathbb{C}[\tau,\alpha,\beta]_{d_\Sigma^{(1)}})\times\mathbb{P}(\mathbb{C}[\tau,\alpha,\beta]_{d_\Sigma^{(2)}})$ be the dense open subset parametrizing generic pairs of curves of degrees $d_\Sigma^{(1)}$ and $d_\Sigma^{(2)}$. Consider the morphism
\begin{align*}
\varphi\colon\mathbb{P}(\mathbb{C}[\tau,\alpha,\beta]_{d_\Sigma^{(1)}})\times\mathbb{P}(\mathbb{C}[\tau,\alpha,\beta]_{d_\Sigma^{(2)}})&\rightarrow\mathbb{P}(\mathbb{C}[\tau,\alpha,\beta]_{d_\Sigma})\\
([f],[g])&\rightarrow[f\cdot g].
\end{align*}
Then, $\varphi(U)/\mathrm{Aut}(\mathbb{P}(1,p,q))$ is the moduli space parametrizing the isomorphism classes of pairs $\left(Y,T+\frac{1}{2}B_0'|_Y\right)$. Moreover, $\varphi$ is injective for $\Sigma\neq W_{15}$, and $2$-to-$1$ if $\Sigma=W_{15}$. In particular,
\[
\dim(\varphi(U))=\dim(U)=\dim(\mathbb{P}(\mathbb{C}[\tau,\alpha,\beta]_{d_\Sigma^{(1)}}))+\dim(\mathbb{P}(\mathbb{C}[\tau,\alpha,\beta]_{d_\Sigma^{(2)}})).
\]
In conclusion, the moduli count for the isomorphism classes of pairs $\left(Y,T+\frac{1}{2}B_0'|_Y\right)$ is
\[
\dim(\mathbb{P}(\mathbb{C}[\tau,\alpha,\beta]_{d_\Sigma^{(1)}}))+\dim(\mathbb{P}(\mathbb{C}[\tau,\alpha,\beta]_{d_\Sigma^{(2)}}))-\dim(\mathrm{Aut}(\mathbb{P}(1,p,q))).
\]
This can be computed in each case $\Sigma\in\{Z_{12},W_{15}\}$, yielding the values in the statement (for the dimension of $\mathrm{Aut}(\mathbb{P}(1,p,q))$ see \cite[Proposition~5.1]{Has00}).

For $\Sigma=Z_{15}$, the moduli count for $\left(Y,T+\frac{1}{2}B_0'|_Y\right)$ is more subtle. First of all, from the construction in Definition~\ref{def:weighted-blow-up-of-family-Sigma} we have that $T=V(\tau)$. However, we should instead think of $T$ as a generic line $\tau+\alpha=0$ passing through the singular point $[0:0:1]$ (in the previous cases, $T$ is unique). Let us write $B_0'|_Y=B_1+B_2$ with $\deg(B_1)=3$ and $\deg(B_2)=4$. As we see from Definition~\ref{def:weighted-blow-up-of-family-Sigma}, the equation defining $B_2$ has all degree $4$ monomials, except $\alpha^4$, which is equivalent to $B_2$ passing through the point $p:=[0:1:0]$. Up to the action of $\mathrm{Aut}(\mathbb{P}(1,1,2))$, we can always assume that $p\in T\cap B_2$, and hence $T$ and $B_2$ are as in Definition~\ref{def:weighted-blow-up-of-family-Sigma}. Therefore, the moduli count for $\left(Y,T+\frac{1}{2}B_0'|_Y\right)$ is given by
\begin{align*}
\dim(&\mathbb{P}(\mathbb{C}[\tau,\alpha,\beta]_1))+\dim(\mathbb{P}(\mathbb{C}[\tau,\alpha,\beta]_3))+\dim(\mathbb{P}(\mathbb{C}[\tau,\alpha,\beta]_4))-\dim(\mathrm{Aut}(\mathbb{P}(1,1,2)))=\\
=(&2-1)+(6-1)+(9-1)-7=7.\qedhere
\end{align*}
\end{proof}

\subsection{Study of the stability of the pair \texorpdfstring{$\left(X_0',\frac{1}{2}B_0'\right)$}{(X_0',1/2B_0')}}
\label{sec:KSBA-stability-of-central-fiber-and-contraction}

\begin{proposition}
\label{prop:central-fiber-after-weighted-blow-up}
Let $\Sigma\in\{Z_{12},W_{15},Z_{15},N_{16}\}$. Then, the pair $\left(X_0',\frac{1}{2}B_0'\right)$ as in Definition~\ref{def:weighted-blow-up-of-family-Sigma} is semi-log canonical. Moreover, the following hold:

\begin{enumerate}

\item If $\Sigma=Z_{12}$, then the pair is stable;

\item If $\Sigma=W_{15},Z_{15}$, then $K_Y+T+\frac{1}{2}B_0'|_Y$ is ample and $K_Z+E+\frac{1}{2}B_0'|_Z$ is nef and big.

\item If $\Sigma=N_{16}$, then $K_Y+T+\frac{1}{2}B_0'|_Y$ is ample and $K_Z+E+\frac{1}{2}B_0'|_Z$ is nef with square $0$.

\end{enumerate}
\end{proposition}

\begin{proof}
The proof structure and techniques are analogous to the ones of \cite[Theorem~3.14]{GPSZ24}, with the needed adjustments as we work with a different type of algebraic surface. We start by proving that $\left(X_0',\frac{1}{2}B_0'\right)$ is semi-log canonical. Also in this case, the key property behind this is that we took a weighted blow up of an isolated semi-quasihomogeneous singularity.

We have that $X_0'$ is a stable toric surface and $-K_{X_0'}$ equals the outside toric boundary \cite[Lemma~2.6.7~(1)]{Ale15}. Hence, $K_{X_0'}$ is $\mathbb{Q}$-Cartier. By construction, $\mathcal{X}'$ has only cyclic quotient singularities, so every Weil divisor is $\mathbb{Q}$-Cartier. In particular $\mathcal{B}'$ is $\mathbb{Q}$-Cartier, which implies that also $B_0'$ is $\mathbb{Q}$-Cartier. In conclusion, $K_{X_0'}+\frac{1}{2}B_0'$ is $\mathbb{Q}$-Cartier.

In $Z\setminus E\cong\mathbb{P}^2\setminus\{\xi\}$ the divisor $B_0'|_Z$ is simple normal crossing by Lemma~\ref{lem:Sigma-singularity-at-the-limit}. So, the pair $\left(Z,E+\frac{1}{2}B_0'|_Z\right)$ log canonical away from $E$, and it is log canonical in a neighborhood of $E$ by an analogous argument as in the proof of \cite[Proposition~3.22]{GPSZ24}. In conclusion, $\left(Z,E+\frac{1}{2}B_0'|_Z\right)$ is log canonical.

The pair $\left(Y,T+\frac{1}{2}B_0'|_Y\right)$ is log canonical away from $T$ because in this open subset the divisor $B_0'|_Y$ is simple normal crossing by Lemma~\ref{lem:moduli-for-curves-in-P1pq}~(1). Moreover, the intersection $B_0'|_Y\cap T$ consists of reduced points, as we verified case by case in Definition~\ref{def:weighted-blow-up-of-family-Sigma}. Then, $\left(Y,T+\frac{1}{2}B_0'|_Y\right)$ is also log canonical in a neighborhood of $T$ (see the conclusion is as in the proof of \cite[Proposition~3.21]{GPSZ24}). This shows that $\left(X_0',\frac{1}{2}B_0'\right)$ is semi-log canonical.

To have that $\left(X_0',\frac{1}{2}B_0'\right)$ is a stable pair, we have to check the ampleness of $K_{X_0'}+\frac{1}{2}B_0'$. This is equivalent to check the ampleness of
\[
K_{Y}+T+\frac{1}{2}B_0'|_Y~\textrm{and}~K_{Z}+E+\frac{1}{2}B_0'|_Z.
\]
The former is ample as
\[
K_Y+T+\frac{1}{2}B_0'|_Y\sim_\mathbb{Q}(-T-D_\alpha-D_\beta)+T+\frac{1}{2}d_\Sigma T\sim_\mathbb{Q}\left(\frac{1}{2}d_\Sigma-p-q\right)T,
\]
where $D_\alpha=V(\alpha),D_\beta=V(\beta)$ and $\frac{1}{2}d_\Sigma-p-q>0$ in our cases. To study the ampleness of the latter, let $D_x=V(x)$, $D_y=V(y)$, $D_z=V(z)$ be the torus fixed curves of $\mathbb{P}^2$ and denote by $\widehat{D}_x,\widehat{D}_y,\widehat{D}_z$ their respective strict transforms in $Z$. Then, $K_{Z}+E+\frac{1}{2}B_0'|_Z$ is ample if and only if its intersections with $\widehat{D}_x,\widehat{D}_y,\widehat{D}_z,E$ are positive. These intersections are computed as in the proof of \cite[Proposition~3.25]{GPSZ24}, together with the equivalence:
\begin{align*}
K_Z+E+\frac{1}{2}B_0'|_Z&\sim_\mathbb{Q}(-\widehat{D}_x-\widehat{D}_y-\widehat{D}_z-E)+E+\frac{1}{2}(7\widehat{D}_x-d_\Sigma E)\\
&\sim_\mathbb{Q}(-\widehat{D}_x-\widehat{D}_x+pE-\widehat{D}_x+qE-E)+E+\frac{1}{2}(7\widehat{D}_x-d_\Sigma E)\\
&=\frac{1}{2}\widehat{D}_x+\left(p+q-\frac{1}{2}d_\Sigma\right)E.
\end{align*}
While the intersections with $\widehat{D}_x$ and $E$ are positive for all $\Sigma$ in our setting, the intersections with $\widehat{D}_y,\widehat{D}_z$ may be zero. More in detail:
\begin{align*}
\left(K_{Z}+E+\frac{1}{2}B_0'|_Z\right)\cdot\widehat{D}_y&=\frac{1}{2}+\left(p+q-\frac{1}{2}d_\Sigma\right)\cdot\frac{1}{q},\\
\left(K_{Z}+E+\frac{1}{2}B_0'|_Z\right)\cdot\widehat{D}_z&=\frac{1}{2}+\left(p+q-\frac{1}{2}d_\Sigma\right)\cdot\frac{1}{p}.
\end{align*}
These two quantities are positive if $\Sigma=Z_{12}$, from which we can conclude that the pair $\left(X_0',\frac{1}{2}B_0'\right)$ is stable. If $\Sigma=W_{15},Z_{15}$, then the intersection with $\widehat{D}_y$ is positive, while the intersection with $\widehat{D}_z$ is zero. For $\Sigma=N_{16}$, both the above intersection numbers are zero. Therefore, for $\Sigma\neq Z_{12}$, we have that $K_{Z}+E+\frac{1}{2}B_0'|_Z$ is nef. Finally, we have that
\[
\left(K_{Z}+E+\frac{1}{2}B_0'|_Z\right)^2=\frac{1}{4}-\frac{1}{pq}\left(p+q-\frac{1}{2}d_\Sigma\right)^2.
\]
This equals $\frac{1}{12},\frac{1}{8},0$ if $\Sigma$ is respectively $W_{15},Z_{15},N_{16}$. Hence, for $\Sigma=W_{15},Z_{15}$, we have that $K_{Z}+E+\frac{1}{2}B_0'|_Z$ is big as it is nef with positive self-intersection.
\end{proof}

\begin{corollary}
\label{cor:threefolds-klt}
Let $\Sigma\in\{Z_{12},W_{15},Z_{15},N_{16}\}$ and let $\left(\mathcal{X}',\frac{1}{2}\mathcal{B}'\right)\rightarrow\Delta$ be a one-parameter family as in Definition~\ref{def:weighted-blow-up-of-family-Sigma}. Also, let $\mathcal{S}'\rightarrow\mathcal{X}'$ be the $(\mathbb{Z}/2\mathbb{Z})^2$-cover with branch divisor $\mathcal{B}'$ with Kunev surfaces as general fibers. Then, the threefold pairs $\left(\mathcal{X}',\frac{1}{2}\mathcal{B}'\right)$ and $\left(\mathcal{S}',0\right)$ are klt. Moreover, the one parameter family $(\mathcal{S}',0)\rightarrow\Delta$ is the semi-stable minimal model in the sense of \cite[Theorem~7.9]{KM98}.
\end{corollary}

\begin{proof}
We prove that $\left(\mathcal{X}',\frac{1}{2}\mathcal{B}'\right)$ is klt by \cite[Lemma~7.2~(5)]{KM98}. To apply it, we first need to check that the morphism
\[
f\colon\left(\mathcal{X}',\frac{1}{2}\mathcal{B}'\right)\rightarrow\Delta
\]
is log canonical: see \cite[Definition~7.1]{KM98}. So, let $X_t'$ be any fiber of $f$. For $t=0$, we have that $\left(X_0',\frac{1}{2}B_0'\right)$ is semi-log canonical by Proposition~\ref{prop:central-fiber-after-weighted-blow-up}. For $t\neq0$ we also have that $\left(X_t',\frac{1}{2}B_t'\right)$ is log canonical (actually more: klt) because $X_t'=\mathbb{P}^2$ and $B_t'=C_1+C_2+L$ is simple normal crossing. Then, by the inversion of adjunction \cite{Kaw07}, we have that $\left(\mathcal{X}',X_t'+\frac{1}{2}\mathcal{B}'\right)$ is log canonical in a neighborhood of $X_t'$ for all $t$. This also implies that for all $t$, the pair $\left(\mathcal{X}',X_t'+\frac{1}{2}\mathcal{B}'\right)$ is log canonical, from which we can conclude that the morphism $f$ is log canonical. Moreover, since for all $t\neq0$ the fiber $\left(X_t',\frac{1}{2}B_t'\right)$ is klt, then $\left(\mathcal{X}',\frac{1}{2}\mathcal{B}'\right)$ is klt by \cite[Lemma~7.2~(5)]{KM98}. The fact that $\left(\mathcal{X}',\frac{1}{2}\mathcal{B}'\right)$ is klt also implies that $\left(\mathcal{S}',0\right)$ is klt by \cite[Proposition~5.20~(4)]{KM98}.

Finally, to prove that $g\colon(\mathcal{S}',0)\rightarrow\Delta$ is the semi-stable minimal model in the sense of \cite[Theorem~7.9]{KM98}, the main thing to verify is that the morphism $g$ is dlt, again, according to \cite[Definition~7.1]{KM98}. So, for all $t$, we need to show that $(\mathcal{S}',S_t')$ is dlt. For $t\neq0$, $S_t'$ is smooth and $\mathcal{S}'$ is smooth away from the central fiber. So, the only thing we need to verify is that, for $t=0$, the pair $(\mathcal{S}',S_0')$ is dlt in a neighborhood of $S_0'$. This comes from the fact that $S_0'$ has two components with a smooth gluing curve and that $S_0'$ has semi-log-terminal singularities (cf. \cite[Theorem~4.23]{KSB88}) by our explicit weighted blow-up construction.
\end{proof}

For $\Sigma=W_{15},Z_{15},N_{16}$, the central fiber $\left(X_0',\frac{1}{2}B_0'\right)$ of $\left(\mathcal{X}',\frac{1}{2}\mathcal{B}'\right)\rightarrow\Delta$ is not stable and we need to contract appropriate curves in $X_0'$ to obtain the stable model. This is carried out in the next proposition for $W_{15}$ and $Z_{15}$. As $N_{16}$ will exhibit a sensibly different behavior, this is described separately in Section~\ref{sec:stable-non-normal-K3-limits-from-N16}.

\begin{proposition}
\label{prop:contraction-W15-Z1,0-case}
For $\Sigma\in\{W_{15},Z_{15}\}$, let $f\colon\left(\mathcal{X}',\frac{1}{2}\mathcal{B}'\right)\rightarrow\Delta$ be the one-parameter family constructed in Definition~\ref{def:weighted-blow-up-of-family-Sigma}. Then, there exists a birational morphism $\mathcal{X}'\rightarrow\mathcal{X}''$ which only contracts the curve $\widehat{D}_z$ in the component $Z$ of the central fiber. In particular, the central fiber $\left(X_0'',\frac{1}{2}B_0''\right)$ of $\left(\mathcal{X}'',\frac{1}{2}\mathcal{B}''\right)\rightarrow\Delta$ is a stable pair, where $\mathcal{B}''$ is the image of $\mathcal{B}'$ under the contraction. $\left(X_0'',\frac{1}{2}B_0''\right)$ is given by the gluing of $\left(Y,T+\frac{1}{2}B_0''|_Y\right)$ and $\left(Z_c,E+\frac{1}{2}B_0''|_{Z_c}\right)$ along $T\cong\mathbb{P}^1\cong E$, where $Z\rightarrow Z_c$ is the contraction of the curve $\widehat{D}_z$.
\end{proposition}

\begin{proof}
This will be a consequence of the Basepoint-free Theorem \cite[Theorem~3.3]{KM98}. To apply it, we first observe that $\left(\mathcal{X}',\frac{1}{2}\mathcal{B}'\right)$ is klt by Corollary~\ref{cor:threefolds-klt}. Next, we observe that the log canonical divisor $K_{\mathcal{X}'}+\frac{1}{2}\mathcal{B}'$ restricted to $X_t'$ is ample for $t\neq0$ and nef and big for $t=0$ by Proposition~\ref{prop:central-fiber-after-weighted-blow-up}. Since $\Delta$ is the spectrum of a discrete valuation ring, we have that $K_{\mathcal{X}'}+\frac{1}{2}\mathcal{B}'$ is nef and big. Therefore, by the Basepoint-free Theorem \cite[Theorem~3.3]{KM98}, we have that the linear system $\left|b\left(K_{\mathcal{X}'}+\frac{1}{2}\mathcal{B}'\right)\right|$ has no basepoints for $b\gg0$, hence inducing a contraction $\mathcal{X}'\rightarrow\mathcal{X}''$ with exceptional locus given by the curve $\widehat{D}_z$.
\end{proof}

\subsection{Description of the $(\mathbb{Z}/2\mathbb{Z})^2$-cover of $Y$}
\label{sec:description-of-the-cover-of-Y}

The goal of this section is to describe the $(\mathbb{Z}/2\mathbb{Z})^2$-cover $\widetilde{Y}$ of the pair $\left(Y,T+\frac{1}{2}B_0'|_Y\right)$, which is stable for all $\Sigma\in\{Z_{12},W_{15},Z_{15},N_{16}\}$. We have that $\widetilde{Y}$ is the actual irreducible exceptional surface in the stable limit of Kunev surfaces. We analyze each case individually. In each case, the conclusion is that $\widetilde{Y}$ is an ADE K3 surface which is a codimension $2$ complete intersection in weighted projective $4$-dimensional space.

($\Sigma=Z_{12}$) We have that $Y\cong\mathbb{P}(1,2,3)$, and $B_0'|_Y=V(f\cdot g)$, where $f(\tau,\alpha,\beta)$ is a degree five polynomial and $g(\tau,\alpha,\beta)$ a degree six polynomial. To construct the $(\mathbb{Z}/2\mathbb{Z})^2$-cover $\widetilde{Y}$, we use the isomorphism
\begin{align*}
\mathbb{P}(1,2,3)&\rightarrow\mathbb{P}(1,4,6)\\
[\tau:\alpha:\beta]&\mapsto[\overline{\tau}:\overline{\alpha}:\overline{\beta}]:=[\tau:\alpha^2:\beta^2],
\end{align*}
under which the polynomials $f$ and $g$ become of degrees $10$ and $12$, respectively. Then, $\widetilde{Y}$ can then be viewed as the complete intersection of the following two hypersurfaces in $\mathbb{P}(1,4,6,5,6)$, with coordinates $[\overline{\tau}:\overline{\alpha}:\overline{\beta}:u:v]$:
\begin{displaymath}
\widetilde{Y}\colon\left\{ \begin{array}{ll}
u^2=f(\overline{\tau}^2,\overline{\alpha},\overline{\beta})\\
v^2=g(\overline{\tau}^2,\overline{\alpha},\overline{\beta}).
\end{array} \right.
\end{displaymath}
This is the codimension $2$ ADE K3 surface in \cite[Table~2, \#42]{IF00}. We now study the gluing curve $\widetilde{T}\subseteq\widetilde{Y}$, which is the preimage of $T=V(\overline{\tau})\subseteq\mathbb{P}(1,4,6)$ under the covering map. The curve $\widetilde{T}$ is described by the solutions of the following system of equations in $\mathbb{P}(4,6,5,6)$:
\begin{displaymath}
\widetilde{T}\colon\left\{ \begin{array}{ll}
u^2=f(0,\overline{\alpha},\overline{\beta})=\overline{\alpha}\overline{\beta}\\
v^2=g(0,\overline{\alpha},\overline{\beta})=\overline{\alpha}^3+\overline{\beta}^2.
\end{array} \right.
\end{displaymath}
Via the isomorphism $\mathbb{P}(4,6,5,6)\cong\mathbb{P}(2,3,5,3)$, the system is equivalent to
\begin{displaymath}
\widetilde{T}\colon\left\{ \begin{array}{ll}
u=\overline{\alpha}\overline{\beta}\\
v^2=\overline{\alpha}^3+\overline{\beta}^2.
\end{array} \right.
\end{displaymath}
(For simplicity of notation, the coordinates in $\mathbb{P}(2,3,5,3)$ are denoted as for $\mathbb{P}(4,6,5,6)$.) The solutions are then
\[
\left[\overline{\alpha}:\overline{\beta}:\overline{\alpha}\overline{\beta}:\pm\sqrt{\overline{\alpha}^3+\overline{\beta}^2}\right].
\]
So, $\widetilde{T}\rightarrow T$ is of degree $2$ and the branch locus is given by the unique solution to $\overline{\alpha}^3+\overline{\beta}^2=0$ in $\mathbb{P}(2,3)$ and the point $[\overline{\alpha}:\overline{\beta}]=[1:0]$. So, $\widetilde{T}\cong\mathbb{P}^1$ by the Hurwitz formula.

($\Sigma=W_{15}$) The $(\mathbb{Z}/2\mathbb{Z})^2$-cover can then be viewed as the complete intersection of the following two general weighted degree $6$ hypersurfaces in $\mathbb{P}(1,2,3,3,3)$, with coordinates $[\tau:\alpha:\beta:u:v]$:
\begin{displaymath}
\widetilde{Y}\colon\left\{ \begin{array}{ll}
u^2=f(\tau,\alpha,\beta)\\
v^2=g(\tau,\alpha,\beta).
\end{array} \right.
\end{displaymath}
Then, $\widetilde{Y}$ is the codimension $2$ ADE K3 surface in \cite[Table~2, \#13]{IF00}. Let us study the preimage $\widetilde{T}$ of $T$ under the cover $\widetilde{Y}\rightarrow Y$. We have that $\widetilde{T}$ is given by the solutions of the system
\begin{displaymath}
\widetilde{T}\colon\left\{ \begin{array}{ll}
u^2=f(0,\alpha,\beta)=a_1\alpha^3+b_1\beta^2\\
v^2=g(0,\alpha,\beta)=a_2\alpha^3+b_2\beta^2.
\end{array} \right.
\end{displaymath}
(We use the coefficients to distinguish the two equations.) The solutions are then
\[
\left[\alpha:\beta:\pm\sqrt{a_1\alpha^3+b_1\beta^2}:\pm\sqrt{a_2\alpha^3+b_2\beta^2}\right].
\]
So, the cover $\widetilde{T}\rightarrow T$ has exactly two preimages over the unique solution to $b_1\beta^2+a_1\alpha^3=0$ (resp. $b_2\beta^2+a_2\alpha^3=0$) and over $[\alpha:\beta]=[1:0]$. Outside of these three points, $\widetilde{T}\rightarrow T$ is $4$-to-$1$. By the Hurwitz formula we can conclude that $\widetilde{T}\cong\mathbb{P}^1$.

($\Sigma=Z_{15}$) We have that $Y\cong\mathbb{P}(1,1,2)$, and $B_0'|_Y=V(f\cdot g)$, where $f(\tau,\alpha,\beta)$ is a degree four polynomial that is general among those without the monomial $\alpha^4$ and $g(\tau,\alpha,\beta)$ a general degree three polynomial. To describe the $(\mathbb{Z}/2\mathbb{Z})^2$-cover, we use the isomorphism
\begin{align*}
\mathbb{P}(1,1,2)&\rightarrow\mathbb{P}(1,2,4)\\
[\tau:\alpha:\beta]&\mapsto[\overline{\tau}:\overline{\alpha}:\overline{\beta}]:=[\tau:\alpha^2:\beta^2],
\end{align*}
under which the polynomials $f$ and $g$ become of degrees $8$ and $6$, respectively. Then, $\widetilde{Y}$ can then be viewed as the complete intersection of the following two hypersurfaces in $\mathbb{P}(1,2,4,4,3)$, with coordinates $[\overline{\tau}:\overline{\alpha}:\overline{\beta}:u:v]$:
\begin{displaymath}
\widetilde{Y}\colon\left\{ \begin{array}{ll}
u^2=f(\overline{\tau}^2,\overline{\alpha},\overline{\beta})\\
v^2=g(\overline{\tau}^2,\overline{\alpha},\overline{\beta}).
\end{array} \right.
\end{displaymath}
This is the codimension $2$ ADE K3 surface in \cite[Table~2, \#19]{IF00}. For the preimage $\widetilde{T}$ of $T$ under the cover $\widetilde{Y}\rightarrow Y$, We have that $\widetilde{T}$ is given by the solutions of the following system of equations in $\mathbb{P}(2,4,4,3)$:
\begin{displaymath}
\widetilde{T}\colon\left\{ \begin{array}{ll}
u^2=f(0,\overline{\alpha},\overline{\beta})=a_1\overline{\alpha}^2\overline{\beta}+b_1\overline{\beta}^2\\
v^2=g(0,\overline{\alpha},\overline{\beta})=a_2\overline{\alpha}^3+b_2\overline{\alpha}\overline{\beta}.
\end{array} \right.
\end{displaymath}
Via the isomorphism $\mathbb{P}(2,4,4,3)\cong\mathbb{P}(1,2,2,3)$ the solutions correspond to the solutions of
\begin{displaymath}
\widetilde{Y}\colon\left\{\begin{array}{ll}
u^2=\overline{\beta}(b_1\overline{\beta}+a_1\overline{\alpha}^2)\\
v=\overline{\alpha}(b_2\overline{\beta}+a_2\overline{\alpha}^2).
\end{array} \right.
\end{displaymath}
The solutions are then
\[
\left[\overline{\alpha}:\overline{\beta}:\pm\sqrt{\overline{\beta}(b_1\overline{\beta}+a_1\overline{\alpha}^2)}:\overline{\alpha}(b_2\overline{\beta}+a_2\overline{\alpha}^2)\right].
\]
So, generically over $T$ we have two preimages, and we have a unique preimage over the points
\[
[1:0:0:a_2],~[\sqrt{b_1}:-a_1:0:\sqrt{b_1}(a_2b_1-a_1b_2)]=[-\sqrt{b_1}:-a_1:0:-\sqrt{b_1}(a_2b_1-a_1b_2)].
\]
Therefore, also in this case the preimage of $T$ under the $(\mathbb{Z}/2\mathbb{Z})^2$-cover is isomorphic to $\mathbb{P}^1$.

($\Sigma=N_{16}$) we have that $Y\cong\mathbb{P}(1,1,1)=\mathbb{P}^2$, and $B_0'|_Y=V(f\cdot g)$, where $f(\tau,\alpha,\beta)$ is a degree two polynomial and $g(\tau,\alpha,\beta)$ a degree three polynomial. To describe the $(\mathbb{Z}/2\mathbb{Z})^2$-cover, we use the isomorphism $\mathbb{P}^2\cong\mathbb{P}(1,2,2)$ under which the polynomials $f$ and $g$ become of degrees $4$ and $6$, respectively. Then, $\widetilde{Y}$ can then be viewed as the complete intersection of the following two hypersurfaces in $\mathbb{P}(1,2,2,2,3)$, with coordinates $[\overline{\tau}:\overline{\alpha}:\overline{\beta}:u:v]$:
\begin{displaymath}
\widetilde{Y}\colon\left\{ \begin{array}{ll}
u^2=f(\overline{\tau}^2,\overline{\alpha},\overline{\beta})\\
v^2=g(\overline{\tau}^2,\overline{\alpha},\overline{\beta}).
\end{array} \right.
\end{displaymath}
This is the codimension $2$ ADE K3 surface in \cite[Table~2, \#8]{IF00}. For the preimage $\widetilde{T}$ of $T$ under the cover $\widetilde{Y}\rightarrow Y$, We have that $\widetilde{T}$ is given by the solutions of the following system of equations in $\mathbb{P}(2,2,2,3)$:
\begin{displaymath}
\widetilde{T}\colon\left\{ \begin{array}{ll}
u^2=f(0,\overline{\alpha},\overline{\beta})=\overline{\alpha}^2+\overline{\alpha}\overline{\beta}+\overline{\beta}^2\\
v^2=g(0,\overline{\alpha},\overline{\beta})=\overline{\alpha}^3+\overline{\alpha}^2\overline{\beta}+\overline{\alpha}\overline{\beta}^2+\overline{\beta}^3.
\end{array} \right.
\end{displaymath}
Via the isomorphism $\mathbb{P}(2,2,2,3)\cong\mathbb{P}(1,1,1,3)$ the solutions correspond to the solutions of
\begin{displaymath}
\widetilde{Y}\colon\left\{\begin{array}{ll}
u^2=\overline{\alpha}^2+\overline{\alpha}\overline{\beta}+\overline{\beta}^2\\
v=\overline{\alpha}^3+\overline{\alpha}^2\overline{\beta}+\overline{\alpha}\overline{\beta}^2+\overline{\beta}^3.
\end{array} \right.
\end{displaymath}
So, generically over $T$, the cover $\widetilde{T}\rightarrow T$ is $2:1$ and also in this case the preimage of $T$ under the $(\mathbb{Z}/2\mathbb{Z})^2$-cover is isomorphic to $\mathbb{P}^1$.

\section{The boundary divisors associated to \texorpdfstring{$Z_{12}$}{Z12}, \texorpdfstring{$W_{15}$}{W15}, \texorpdfstring{$Z_{15}$}{Z15}}
\label{sec:boundary-divisors-associated-to-Sigma}

The following construction is an adaptation of \cite[Section~4]{GPSZ24}. Let $\Delta=\mathrm{Spec}(\mathbb{C}[[t]])$ and let us consider $\Sigma\in\{Z_{12},W_{15},Z_{15}\}$. Recall the corresponding weight functions $\mathrm{wt}_1,\mathrm{wt}_2$ and the $\mathbb{C}^*$-actions $\star_1,\star_2$ in Definition~\ref{def:weight-functions}. Let us view $\mathbb{P}^9$ and $\mathbb{P}^2$ respectively as the parameter spaces for cubics and lines in $\mathbb{P}^2$. Let $U\subseteq\mathbb{P}^9\times\mathbb{P}^9\times\mathbb{P}^2$ be the dense open subset such that, if $u=([F_1],[F_2],[L])\in U$, then
\[
F_u(t,[x:y:z]):=((t\star_1F_1)\cdot(t\star_2F_2)\cdot L)([x:y:z])=0
\]
defines a one-parameter family of plane curves as in Definition~\ref{def:1-parameter-degeneration-Sigma} which for $t=0$ is a $\Sigma$ arrangement with $\Sigma$ singularity at $[1:0:0]$. Let
\[
\mathfrak{F}:=U\times\Delta\times\mathbb{P}^2,
\]
which we view as a family of $\mathbb{P}^2$ over $U\times\Delta$. Let $\mathfrak{D}\subseteq\mathfrak{F}$ be the divisor given by
\[
\mathfrak{D}:=V(F_u(t,[x:y:z]))\subseteq\mathfrak{F}.
\]
Therefore, the fiber of $\pi\colon\left(\mathfrak{F},\frac{1}{2}\mathfrak{D}\right)\rightarrow U\times\Delta$ over $(u,t)\in U\times\Delta$ is the pair $\left(\mathbb{P}^2,\frac{1}{2}\mathfrak{D}_{(u,t)}\right)$, and if we specialize for $t\to0$ the plane curve $\mathfrak{D}_{(u,0)}$ has a $\Sigma$ singularity at the point $[1:0:0]$. Both the divisors $\mathfrak{D}$ and $K_\mathfrak{F}$ are Cartier.

In the open patch where $x\neq0$, let $\mathfrak{C}$ be the subscheme of $\mathfrak{F}$ associated to the ideal $(t^{d_\Sigma},y^q,z^p)$. Then, we define $\mathfrak{F}':=\mathrm{Bl}_\mathfrak{C}\mathfrak{F}$, and let $\mathfrak{D}'\subseteq\mathfrak{F}'$ be the strict transform of $\mathfrak{D}$. Then, the fiber of $\pi'\colon\left(\mathfrak{F}',\frac{1}{2}\mathfrak{D}'\right)\rightarrow U\times\Delta$ over $(u,t)$ is a stable pair, and in the special case where $t=0$ the fiber over $(u,0)$ is the gluing of $\left(Y,T+\frac{1}{2}B_0'|_Y\right)$ and $\left(Z,E+\frac{1}{2}B_0'|_Z\right)$ along $T\cong\mathbb{P}^1\cong E$ as in Definition~\ref{def:weighted-blow-up-of-family-Sigma}. We have that $K_{\mathfrak{F}'}$ and $\mathfrak{D}'$ are both $\mathbb{Q}$-Cartier.

In particular, for $N$ large enough and divisible by $p,q,d_\Sigma$, we have that $N\left(K_{\mathfrak{F}'}+\frac{1}{2}\mathfrak{D}'\right)$ is Cartier and it restricts to the fibers $F$ of the family $\mathfrak{F}'\rightarrow U\times\Delta$ giving the Cartier divisor $N\left(K_F+\frac{1}{2}\mathfrak{D}'|_F\right)$. Both $\pi'$ and $\pi'|_{\mathfrak{D}'}$ are flat as they are dominant morphisms from integral schemes to normal schemes with reduced fibers of constant dimension \cite[Lemma~10.12]{HKT09}.

If $\Sigma=Z_{12}$, then the fibers of $\pi'$ are stable pairs: this is clear for the fibers over $(u,t)$ with $t\neq0$, while for the fibers over $(u,0)$ follows from Proposition~\ref{prop:central-fiber-after-weighted-blow-up}~(1). In particular, $\pi'$ is a well-defined family of stable pairs, and therefore it induces a morphism $f_{Z_{12}}\colon U\times\Delta\rightarrow\overline{\mathrm{M}}_{\mathrm{K}}$ to the KSBA compactification of the moduli space of Kunev surfaces.

On the other hand, if $\Sigma=W_{15},Z_{15}$, then the fibers over the points of the form $(u,0)$ are semi-log canonical, $K_Y+T+\frac{1}{2}B_0'|_T$ is ample, but $K_Z+E+\frac{1}{2}B_0'|_Z$ is only nef and big (see Proposition~\ref{prop:central-fiber-after-weighted-blow-up}~(2)), making the pair not stable. Recall that the only obstruction to the ampleness of $K_Z+E+\frac{1}{2}B_0'|_Z$ is the curve $\widehat{D}_z\subseteq Z$ which intersects $K_Z+E+\frac{1}{2}B_0'|_Z$ by $0$. However, by the same argument of Proposition~\ref{prop:contraction-W15-Z1,0-case} we have by the basepoint-free theorem a birational contraction $\mathfrak{F}'\rightarrow\mathfrak{F}''$ that contracts precisely the curves $\widehat{D}_z$ in the fibers over points of the form $(u,0)$. Denote by $\mathfrak{D}''$ the image of $\mathfrak{D}'$ under this morphism. We obtain therefore a commutative diagram:
\begin{center}
\begin{tikzpicture}[>=angle 90]
\matrix(a)[matrix of math nodes,
row sep=2em, column sep=2em,
text height=1.5ex, text depth=0.25ex]
{\left(\mathfrak{F}',\frac{1}{2}\mathfrak{D}'\right)&&\left(\mathfrak{F}'',\frac{1}{2}\mathfrak{D}''\right)\\
&U\times\Delta,&\\};
\path[->] (a-1-1) edge node[above]{}(a-1-3);
\path[->] (a-1-1) edge node[below left]{$\pi'$}(a-2-2);
\path[->] (a-1-3) edge node[below right]{$\pi''$}(a-2-2);
\end{tikzpicture}
\end{center}
where now the fibers of $\pi''$ are stable pairs. By the same argument as for $Z_{12}$ we have that $\pi''$ and $\pi''|_{\mathfrak{D}''}$ are flat and $N\left(K_{\mathfrak{F}''}+\frac{1}{2}\mathfrak{D}''\right)$ is Cartier for $N$ large and divisible enough. In particular, for $\Sigma=W_{15},Z_{1,0}$, the family $\pi''$ induces a morphism $f_\Sigma\colon U\times\Delta\rightarrow\overline{\mathrm{M}}_{\mathrm{K}}$ to the KSBA compactification of the moduli space of Kunev surfaces.

\begin{definition}
\label{def:boundary-stratum-associated-to-Z122}
With the notations introduced above, we define the boundary stratum $D(\Sigma)\subseteq\overline{\mathrm{M}}_{\mathrm{K}}$ as the Zariski closure, with the reduced scheme structure, of $f_\Sigma(U\times\{0\})$.
\end{definition}

\begin{theorem}
\label{thm:divisor-for-Z_12^2}
The boundary strata $D(Z_{12}),D(W_{15}),D(Z_{15})\subseteq\overline{\mathrm{M}}_{\mathrm{K}}$ in Definition~\ref{def:boundary-stratum-associated-to-Z122} are divisors.
\end{theorem}

\begin{proof}
The generic surface parametrized by $D(Z_{12})$ is given by the gluing of two pairs $\left(Y,T+\frac{1}{2}B_0'|_Y\right)$ and $\left(Z,E+\frac{1}{2}B_0'|_Z\right)$ as in Definition~\ref{def:weighted-blow-up-of-family-Sigma}.

The moduli count for $\left(Z,E+\frac{1}{2}B_0'|_Z\right)$ is given by
\[
(6-1)+(5-1)+(3-1)-5=6,
\]
where
\begin{itemize}

\item $6$ (resp. $5$) is the number of monomials in $\lim\limits_{t\to0}t\star_1F_1$ (resp. $\lim\limits_{t\to0}t\star_2F_2$), minus one for projective scaling. These are computed in the proof of Lemma~\ref{lem:Sigma-singularity-at-the-limit}. $3-1$ is the parameter count for a generic line in $\mathbb{P}^2$.

\item $5$ is the dimension of the subgroup of $\mathrm{PGL}_3$ fixing the point $\xi=[1:0:0]$ and preserving the monomials $x^iy^jz^k$ such that $\mathrm{wt}_1(x^iy^jz^k)\geq0$ and $\mathrm{wt}_2(x^iy^jz^k)\geq0$.

\end{itemize}
On the other hand, the pair $\left(Y,T+\frac{1}{2}B_0'|_Y\right)$ depends on $5$ parameters up to automorphism by Lemma~\ref{lem:moduli-for-curves-in-P1pq}~(2). As for the gluing of $Y$ and $Z$ along $T\cong\mathbb{P}^1\cong E$, the three marked points $T\cap B_0'|_Y$ and $E\cap B_0'|_Z$ are identified, so there is no contribution to the parameter count coming from the gluing locus. In conclusion,
\[
\dim(D(Z_{12}))=6+5-0=11.
\]

For $W_{15}$ (resp. $Z_{15}$) the argument is mostly analogous to the one for the singularity $Z_{12}$, with a few adaptations. The generic surface $\left(X_0'',\frac{1}{2}B_0''\right)$ parametrized by $D({W_{15}})$ is given by the gluing of two pairs $\left(Y,T+\frac{1}{2}B_0''|_Y\right)$ and $\left(Z_c,E+\frac{1}{2}B_0''|_Z\right)$ as in Proposition~\ref{prop:contraction-W15-Z1,0-case}. We have that $\left(Z_c,E+\frac{1}{2}B_0''|_Z\right)$ depends on $5$ (resp. $5$) parameters up to automorphisms (the birational transformations that brought from $\mathbb{P}^2$ to $Z_c$ did not alter this count). On the other hand, the pair $\left(Y,T+\frac{1}{2}B_0''|_Y\right)$ depends on $7$ (resp. $7$) parameters up to automorphism by Lemma~\ref{lem:moduli-for-curves-in-P1pq}~(2). The gluing of $Y$ and $Z_c$ along $T\cong\mathbb{P}^1\cong E$ can be understood as follows. For $W_{15}$ the two marked points $T\cap B_0'|_Y$ and $E\cap B_0'|_Z$ need to be identified, but also the two singular points of $Y$ on $T$ need to be identified with the two singular points of $Z_c$ on $E$, so we have to account for $1$ parameter needed to have compatible gluing. The same conclusion hold for $Z_{15}$, with the only difference that $T\cap B_0'|_Y$ and $E\cap B_0'|_Z$ already consist of four distinct points. In conclusion,
\[
\dim(D(W_{15}))=5+7-1=11,~\dim(D(Z_{15}))=5+7-1=11.\qedhere
\]
\end{proof}


\section{The boundary stratum associated to \texorpdfstring{$N_{16}$}{N16}}
\label{sec:stable-non-normal-K3-limits-from-N16}

We now go back to the one-parameter family $(\mathcal{X}',\frac{1}{2}\mathcal{B}')\rightarrow\Delta$ arising from the $N_{16}$ singularity and compute the stable replacement of its central fiber.

\begin{proposition}
\label{prop:stable-replacement-for-N16}
Let $f\colon\left(\mathcal{X}',\frac{1}{2}\mathcal{B}'\right)\rightarrow\Delta$ be the one-parameter family constructed in Definition~\ref{def:weighted-blow-up-of-family-Sigma} for $\Sigma=N_{16}$. Let $\mathcal{X}'\rightarrow\mathcal{X}''$ be the contraction of $Z\cong\mathbb{F}_1$ along its ruling and denote by $\mathcal{B}''$ the pushforward of $\mathcal{B}'$ under this contraction. Then, the central fiber of $\left(\mathcal{X}'',\frac{1}{2}\mathcal{B}''\right)\rightarrow\Delta$ is given by $\left(Y,T+\frac{1}{2}B_0'|_Y\right)$, which is stable.
\end{proposition}

\begin{proof}
From the proof of Proposition~\ref{prop:central-fiber-after-weighted-blow-up} we have that the intersection products of $K_Z+E+\frac{1}{2}B_0'|_Z$ with $\widehat{D}_y$ and $\widehat{D}_z$ are equal to zero. As $Z$ is isomorphic to the Hirzebruch surface and $\widehat{D}_y,\widehat{D}_z$ are two linearly equivalent rulings, we have that $\left(K_Z+E+\frac{1}{2}B_0'|_Z\right)\cdot\ell=0$ for any ruling $\ell\subseteq Z$, and these are the only irreducible curves with intersection product zero with $K_Z+E+\frac{1}{2}B_0'|_Z$. Therefore, the stable model of $(\mathcal{X}',\frac{1}{2}\mathcal{B}')\rightarrow\Delta$ is obtained by contracting $Z$ along its ruling. Let $\rho\colon\mathcal{X}'\rightarrow\mathcal{X}''$ be such contraction and let $\rho_0\colon X_0'\rightarrow Y$ be its restriction to the central fiber. Recall that $B_0'|_Z$ consists of the following three types of components:
\begin{itemize}

\item The strict transform of the line $L$, which is a section of $Z$;

\item The strict transform of $C_1$, which is also a section of $Z$;

\item The strict transform of $C_2$, which consists of three distinct rulings of $Z$.

\end{itemize}
Then, $\rho_{0*}(B_0'|_Z)=2T$, and hence $B_0''=\rho_{0*}(B_0')=2T+B_0'|_Y$. From this we obtain that the central fiber of $\left(\mathcal{X}'',\frac{1}{2}\mathcal{B}''\right)\rightarrow\Delta$ is $\left(Y,T+\frac{1}{2}B_0'|_Y\right)$ as claimed. This is stable as we already know that it is log canonical and $K_Y+T+\frac{1}{2}B_0'|_Y$ is ample by Proposition~\ref{prop:central-fiber-after-weighted-blow-up}.
\end{proof}

\begin{figure}
\centering
\begin{tikzpicture}
\draw [line width=1.5pt,dash pattern=on 2pt off 2pt]  plot[domain=-1.8:1.8,variable=\t]({\t},{\t^3-\t^2-\t});
\draw [rotate around={0.:(0.,0.)},line width=1.pt] (0.,0.) ellipse (2.59615659807322cm and 0.7000207723483153cm);
\draw [line width=1.0pt] (-3.,0.04)-- (3.,0.04);
\draw [line width=1.0pt] (-3.,-0.04)-- (3.,-0.04);

\draw (-0.8,-0.3) node {{\scriptsize$L_1=L$}};
\draw (0.8,1.1) node {$Q_1$};
\draw (1.6,-1.2) node {$C_2$};

\end{tikzpicture}
\caption{Plane curve arrangement $C_1+C_2+L$ arising from the stable replacement in the $N_{16}$ case.}
\label{fig:double-line-arrangement}
\end{figure}
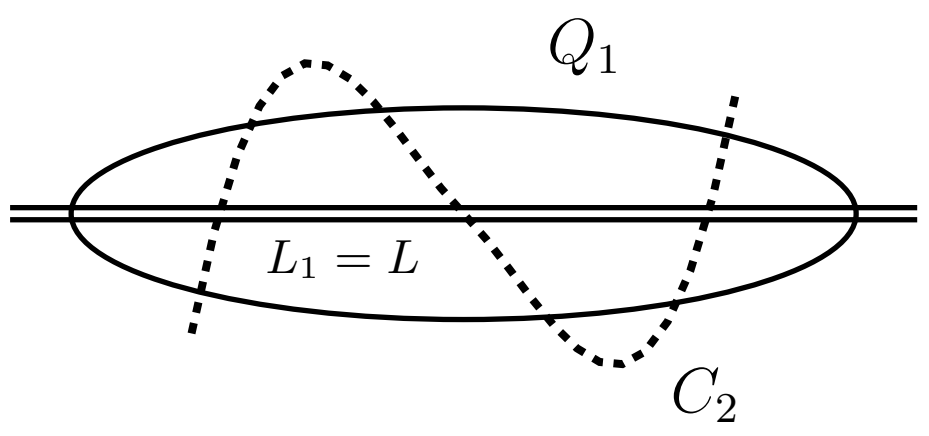

Let us consider the stable central fiber $\left(X_0'',\frac{1}{2}B_0''\right)$ as in Proposition~\ref{prop:stable-replacement-for-N16}. Then, $X_0''\cong\mathbb{P}^2$, and $B_0''$ consists of the union of a smooth conic, a smooth cubic (see the description in Definition~\ref{def:weighted-blow-up-of-family-Sigma}), and a double line. More precisely, we can write $B_0''=C_1+C_2+L$ such that $C_1=Q_1+L$ with $Q_1$ smooth conic and $L$ line, $C_2$ is a smooth cubic, and the support of $C_1+C_2+L$ is normal crossing: see Figure~\ref{fig:double-line-arrangement}. From this we see that the pair $\left(\mathbb{P}^2,\frac{1}{2}(C_1+C_2+L)\right)$ is log canonical. This case already appears in \cite[Example~5.4]{FPR17}.

Next, we formalize the definition of the boundary stratum in $\overline{\mathrm{M}}_{\mathrm{K}}$ generically parametrizing such stable degenerations arising from $N_{16}$.

\begin{definition}
\label{def:boundary-stratum-associated-to-N16}
Denote by $S(N_{16})$ the locally closed subset of
\[
\mathbb{P}(H^0(\mathbb{P}^2,\mathcal{O}(3)))\times\mathbb{P}(H^0(\mathbb{P}^2,\mathcal{O}(3)))\times\mathbb{P}(H^0(\mathbb{P}^2,\mathcal{O}(1)))
\]
parametrizing curve arrangements $C_1+C_2+L$, where one of the two cubics is smooth, the other is the union of a smooth conic and a line $L$, and the support of $C_1+C_2+L$ is normal crossing. Then, we have that the open inclusion $\iota\colon\mathrm{M}_{\mathrm{K}}\hookrightarrow\overline{\mathrm{M}}_{\mathrm{K}}$ constructed in Section~\ref{sec:construction-of-the-compactification} extends to $S(N_{16})/(\mathrm{SL}_3\times S_2)$, where recall $S_2$ acts permuting the two copies of $\mathbb{P}(H^0(\mathbb{P}^2,\mathcal{O}(3)))$. We denote by $D(N_{16})$ the boundary stratum given by the closure of $S(N_{16})/(\mathrm{SL}_3\times S_2)$ in $\overline{\mathrm{M}}_{\mathrm{K}}$.
\end{definition}

The next proposition then follows by counting the moduli of the plane curve $C_1+C_2+L$ as in Definition~\ref{def:boundary-stratum-associated-to-N16} up to $(\mathrm{SL}_3\times S_2)$-action.

\begin{proposition}
\label{prop:dimension-count-for-D(N16)}
Let $D(N_{16})\subseteq\overline{\mathrm{M}}_{\mathrm{K}}$ as in Definition~\ref{def:boundary-stratum-associated-to-N16}. Then, $\dim(D(N_{16}))=8$.
\end{proposition}

\begin{remark}
There are other ways of acquiring a double line in $C_1+C_2+L$, which describe other stable surfaces parametrized by $\overline{\mathrm{M}}_{\mathrm{K}}$. These can be found in \cite[Section~5]{FPR17}.
\end{remark}


\section{Classification of K3 tails in \texorpdfstring{$\overline{\mathrm{M}}_{\mathrm{K}}$}{barMK}}
\label{sec:classification-of-K3-tails-for-MKbar}

We now aim to show that the K3 surfaces that arise as exceptional divisors in the way described in Sections~\ref{sec:non-log-canonical-singularities-Z12-W15-Z15-N16} and \ref{sec:stable-degen-of-Kunev-surfaces} are the only possible ones. We make this formal in the next definition.

\begin{definition}
\label{def:K3-tail-in-moduli-of-Kunev-surfaces}
Let $\mathcal{S}\rightarrow\Delta$ be a one-parameter degeneration of Kunev surfaces that, fiberwisely, is the $(\mathbb{Z}/2\mathbb{Z})^2$-cover of $\left(\mathcal{X}=\Delta\times\mathbb{P}^2,\frac{1}{2}\mathcal{B}\right)\rightarrow\Delta$, where
\begin{itemize}

\item for $t\neq0$, $B_t$ is the simple normal crossing union of two smooth cubics and one line;

\item for $t=0$, $B_0=C_1+C_2+C_3$ has one isolated non-log canonical singular point at $\xi=[1:0:0]$ that is a semi-quasihomogeneous singularity, and it is normal crossing elsewhere.

\end{itemize}
Let $\zeta\in S_0$ be the unique isolated singularity corresponding to $\xi$. Let $\mathcal{X}'\rightarrow\mathcal{X}$ be the weighted blow up at $\xi$ resolving the non-log canonical singularity, and consider the following commutative square:
\begin{center}
\begin{tikzpicture}[>=angle 90]
\matrix(a)[matrix of math nodes,
row sep=2em, column sep=2em,
text height=1.5ex, text depth=0.25ex]
{\mathcal{S}'&\mathcal{S}\\
\mathcal{X}'&\mathcal{X},\\};
\path[->] (a-1-1) edge node[above]{}(a-1-2);
\path[->] (a-1-1) edge node[left]{}(a-2-1);
\path[->] (a-2-1) edge node[below]{}(a-2-2);
\path[->] (a-1-2) edge node[right]{}(a-2-2);
\end{tikzpicture}
\end{center}
where $\mathcal{S}'\rightarrow\mathcal{S}$ is a blow up of $\zeta$. Suppose that:
\begin{itemize}

\item[(i)] $\mathcal{S}'\rightarrow\mathcal{S}$ resolves the non-log canonical isolated singularity and provides a semi-stable minimal model of $\mathcal{S}\rightarrow\Delta$, in the sense of \cite[Theorem~7.9]{KM98};

\item[(ii)] The exceptional surface $\widetilde{Y}\subseteq\mathcal{S}'$ is an ADE K3 surface.
\end{itemize}
Then, we call $\widetilde{Y}$ a \emph{K3 tail} of the KSBA moduli space $\overline{\mathrm{M}}_{\mathrm{K}}$.
\end{definition}

\begin{theorem}
\label{thm:classification-of-K3-tails}
Let $\mathcal{S}\rightarrow\Delta$ be a one-parameter degeneration of Kunev surfaces as in Definition~\ref{def:K3-tail-in-moduli-of-Kunev-surfaces}. Then, the K3 tail $\widetilde{Y}\subseteq\mathcal{S}'$ is a codimension $2$ complete intersection in a $4$-dimensional weighted projective space, and the only possible cases are $\#8$ ($N_{16}$), $\#13$ ($W_{15}$), $\#19$ ($Z_{15}$), and $\#42$ ($Z_{12}$) in \cite[Table~2]{IF00}.
\end{theorem}

\begin{proof}
We note that the assumption (i) implies that the total family $\mathcal{S}'$ is terminal with $K_{\mathcal{S}'/\Delta}$ being a $\mathbb{Q}$-Cartier nef and big divisor, and the central fiber of $\mathcal S'$ is dlt.
Let $G=(\mathbb{Z}/2\mathbb{Z})^2$, let $\widetilde{S}_0$ be the strict transform of the original central fiber, and let $\widetilde{C}:=\widetilde{S}_0\cap \widetilde{Y}$. 

The deck transformations fix the coordinates $t,y,z$ and only change the signs of the covering coordinates, so the valuation $\nu$ defining the blow up is invariant. Therefore, the exceptional divisor $\operatorname{Proj}\!\left(\operatorname{gr}_{\nu}\mathcal{O}_{\mathcal{S}}\right)$ is $(\mathbb{Z}/2\mathbb{Z})^2$-invariant, so the action extends to the exceptional component $\widetilde{Y}$.

We now distinguish whether $L$ passes through $\xi$ or not.

\noindent\textbf{Case 1:} $\xi\in L$. We have the following chain of double covers:
\[
\widetilde Y\longrightarrow P'\longrightarrow P, \qquad P=\mathbb{P}(1,p,q)\ \textrm{or}\ \mathbb{P}(1,2p,2q),
\]
where $P'\to P$ is a double cover branched along two curves $C_n,C_m$ of even weighted degrees $n,m$, and $\widetilde Y\to P'$ has branch induced by $L$.
Since $\widetilde Y$ is an ADE K3 surface, the classification of involutions on K3 surfaces implies that $P'$ is rational.
Let $Z'$ be the double cover of $Z=\mathrm{Bl}^{(p,q)}_\xi\mathbb{P}^2$ branched along the strict transforms of $C_1+C_2$. We then have
\[
\widetilde{S}_0\longrightarrow Z'\longrightarrow Z.
\]
Away from the gluing curve $C'=Z'\cap P'$, the surface $Z'$ has only nodes; along the corresponding curve $\widetilde{C}$, 
the surface $\widetilde{S}_0$ has only quotient singularities.

By assumption (i), $\mathcal{S}'\rightarrow\Delta$ is a semi-stable minimal model and the dual complex $\Gamma$ of the central fiber $S_0'$ satisfies $h^1(|\Gamma|)=h^2(|\Gamma|)=0$. Therefore, after passing to the semi-stable reduction of the central fiber and applying \cite[Page~118]{Mor84}, we can argue the following Hodge-theoretic formula:
\begin{equation}
\label{eq:hodge-theoretic-formula-bidouble-cover}
p_g(\widetilde{S}_0)+p_g(\widetilde Y)+g(\widetilde{C})=p_g(S_t)=1.
\end{equation}
An analogous argument can be found in the proof of \cite[Theorem~6.1]{GPSZ24}. \cite[Proof of Theorem~10.6]{Kaw88} or \cite[Lemma~3.5, Lemma~5.2, Corollary~5.4]{Cor94} provide an algorithm for the semi-stable reduction process. The key point of this process is to produce only rational surfaces in the exceptional divisors and to keep $h^1(|\Gamma|)=h^2(|\Gamma|)=0$ as we have in the central fiber of $\mathcal{S}'$. Since $p_g(\widetilde Y)=1$, we have that $g(\widetilde{C})=0$, which also implies that $g(C')=0$ as we have a covering morphism $\widetilde{C}\rightarrow C'$. By applying the same formula to the intermediate double cover degeneration, whose general fiber is an ADE K3 surface given by the double cover of $\mathbb{P}^2$ branched along two cubics, we obtain that
\[
p_g(Z')+p_g(P')+g(C')=1.
\]
Since $P'$ is rational and $g(C')=0$, we obtain that $p_g(Z')=1$. The pullback of canonical forms under the finite map $\widetilde{S}_0\to Z'$ gives that $p_g(\widetilde{S}_0)\geq p_g(Z')=1$. This contradicts \eqref{eq:hodge-theoretic-formula-bidouble-cover} because $p_g(\widetilde{Y})=1$. Thus, Case~1 is impossible.

\noindent\textbf{Case 2:} $\xi\notin L$. Again, $\widetilde Y\longrightarrow P'\longrightarrow P$ is a $(\mathbb{Z}/2\mathbb{Z})^2$-cover of a weighted projective plane $P=\mathbb{P}(1,p,q)$ or $\mathbb{P}(1,2p,2q)$.
A $(\mathbb{Z}/2\mathbb{Z})^2$-cover is described by effective divisors $D_1,D_2,D_3$ line bundles $L_1,L_2,L_3$, equations $f_i$ for $D_i$, and covering coordinates $u_i$, satisfying
\[2L_i\sim D_j+D_k, \qquad D_k+L_k\sim L_i+L_j, \qquad
u_3f_3=u_1u_2.\]
The data are unique because $\operatorname{Pic}(P)$ has no $2$-torsion \cite[Section~2]{Cat84}.
Take $D_3$ to be the divisor induced by $L$. Since $L$ does not pass through the center $\xi$, it has no component on the exceptional weighted plane. Thus $f_3=1$, and
$u_3=u_1u_2$. The third covering coordinate is therefore redundant. If the other two branch curves have weighted degrees $n$ and $m$, then
\[
\widetilde Y=X_{n,m}=
\left\{\begin{aligned}
u^2&=f_n(\tau,\alpha,\beta),\\
v^2&=g_m(\tau,\alpha,\beta)
\end{aligned}\right\}
\subseteq\mathbb{P}(a_0,a_1,a_2,a_3,a_4),
\]
where $n=2a_3,\ m=2a_4$. Hence, $\widetilde Y$ is a codimension $2$ weighted complete intersection in a $4$-dimensional weighted projective space. Let $\iota$ be the involution on $\widetilde{Y}$ with quotient $P'=\widetilde{Y}/\langle\iota\rangle$. First, since $\widetilde{Y}$ is klt, we have that $P'$ is klt.
In what follows, we show that also $P'$ is an ADE K3 surface.

By applying again the Hodge-theoretic formula to $\widetilde{S}_0\cup_{\widetilde{C}}\widetilde{Y}$ and from $p_g(\widetilde Y)=1$, we obtain that $g(\widetilde{C})=0$. As $\widetilde{S}_0$ is rational by Proposition~\ref{prop:rationality-of-non-lc-limit} and since we have a dominant morphism $\widetilde{S}_0\to Z'$, also $Z'$ is rational, and hence $p_g(Z')=0$. If we apply the Hodge-theoretic formula to $Z'\cup_{C'}P'$, we also have that $p_g(Z')+p_g(P')+g(C')=1$, and hence $p_g(P')+g(C')=1$. If $p_g(P')=0$ and $g(C')=1$, then we have a contradiction because we have a covering morphism $\widetilde{C}\rightarrow C'$, so the rational curve $\widetilde{C}$ surjects onto the genus one curve $C'$. Hence, $p_g(P')=1$ and $g(C')=0$. Since $p_g(P')=1$, the nonzero canonical form $\omega$ on $P'$ pulls back via $\pi\colon\widetilde{Y}\rightarrow P'$ to the nowhere vanishing canonical form $\pi^*\omega$ on $\widetilde{Y}$. Therefore, the quotient map $\pi$ has no branch divisor, and $\iota^*(\pi^*\omega)=\pi^*\omega$ because $\pi\circ\iota=\pi$. Thus $\iota$ is a symplectic involution, and it implies $P'$ is an ADE K3 surface.

It remains to classify the possibilities for $P'$. To do this, we use \cite{KL25}, which classified the semi-quasihomogeneous non-log canonical surface singularities whose weighted resolution gives an ADE K3 surface. More precisely, let $\mathcal{K}\rightarrow\Delta$ be the one-parameter degeneration of ADE K3 surfaces that are the double covers of $\mathbb{P}^2$ branched along two cubics. In particular, we have a fiber-wise double covering morphism $\mathcal{S}\rightarrow\mathcal{K}$. The central fiber $K_0\subseteq\mathcal{K}$ has a semi-quasihomogeneous singularity over $\xi\in C_1\cap C_2$. The ADE K3 surface $P'$ is the exceptional divisor of the weighted resolution $\mathcal{K}'\rightarrow\mathcal{K}$ of this singularity. So, $P'$ has to fall within \cite[Tables~1, 2]{KL25}. In our cases of interest, the semi-quasihomogeneous singularity must have a double covering structure. Among these, we can recognize the ones we already identified in this paper: namely, $Z_{12},W_{15},Z_{15},N_{16}$. We need to show that the other cases do not arise. Among these, we only have to consider the ones that have the structure of double cover. We list them in Table~\ref{tbl:table_kerr_laza}, together with the corresponding values of $p,q,d_\Sigma$ and, consequently, the potential weighted projective space $\mathbb{P}(a_0,a_1,a_2,a_3,a_4)$ containing the ADE K3 surface $\widetilde{Y}=X_{n,m}$ as codimension $2$ complete intersection.

At the same time, in \cite[Table~2]{IF00} we have a complete classification of the codimension $2$ complete intersection ADE K3 surfaces in weighted projective space $\mathbb{P}(a_0,a_1,a_2,a_3,a_4)$ such that $\sum_{i=0}^4a_i\leq100$. For instance, among these we find the known cases:
\[\begin{array}{c|c|c}
\text{singularity}&\widetilde Y&\text{Iano-Fletcher number}\\ \hline
N_{16}&X_{4,6}\subseteq\mathbb{P}(1,2,2,2,3)&\#8\\
W_{15}&X_{6,6}\subseteq\mathbb{P}(1,2,3,3,3)&\#13\\
Z_{15}&X_{6,8}\subseteq\mathbb{P}(1,2,4,4,3)&\#19\\
Z_{12}&X_{10,12}\subseteq\mathbb{P}(1,4,6,5,6)&\#42
\end{array}\]

We have that the cases in Table~\ref{tbl:table_kerr_laza} have to appear within \cite[Table~2]{IF00} as well. For instance, for $W_{12}$, we have $a_0+a_1+a_2+a_3+a_4=1+4+5+d_\Sigma=30<100$, and the same holds for all the other cases. By direct comparison with \cite[Table~2]{IF00}, we can exclude all the cases in Table~\ref{tbl:table_kerr_laza} except $E_{16}$, $Y_7$, $Y_{42}$ because the corresponding weighted projective spaces appear in \cite[Table~2]{IF00} for some choices of $a_3,a_4$. (We point out that $\{a_3,a_4\}=\{n/2,m/2\}$. In particular, $n$ and $m$ have to be even). We now exclude these possibilities.
\begin{itemize}

\item For $Y_7$ and $Y_{42}$, we have from the last column of \cite[Table~2]{KL25} that the gluing curve $\widetilde{C}=\widetilde{S}_0\cap\widetilde{Y}$ has genus one, which contradicts the equality $p_g(\widetilde{S}_0)+p_g(\widetilde{Y})+g(\widetilde{C})=1$. So, the singularities $Y_7$ and $Y_{42}$ do not arise in our setting.

\item For $E_{16}$ we have that
\[
\widetilde{Y}=X_{8,10}\subseteq\mathbb{P}(1,2,6,4,5).
\]
By \cite{Arn76} $E_{16}$ (also denoted by $J_{3,0}$) has normal form $(y^9+by^3z^2+z^3)+cy^7z=0$ with $4b^3+27\neq0$. Hence, it has three smooth analytic branches with pairwise local intersection number $3$ at $\xi$. As this cannot be attained using two plane cubics, we argue that $E_{16}$ is not possible.\qedhere

\end{itemize}
\end{proof}

\begin{table}
\caption{The surface singularities in \cite[Table~1, 2]{KL25} different from $Z_{12},W_{15},Z_{15},N_{16}$ that can be potentially acquired by the central fiber $K_0\subseteq\mathcal{K}\rightarrow\Delta$ in the proof of Theorem~\ref{thm:classification-of-K3-tails}.}
\label{tbl:table_kerr_laza}
\renewcommand{\arraystretch}{2}
\begin{tabular}{|c|c|c|c|c|c|}
\hline
\text{sing.}&$\text{germ}$&$(p,q)$&$d_\Sigma$&$n+m$&$\text{potential w.p.s.}$\\
\hline
$W_{12}$&$y^5+z^4$&$(4,5)$&$20$&$20$&$\mathbb{P}(1,4,5,a_3,a_4)$\\
\hline
$E_{14}$&$y^8+z^3$&$(3,8)$&$24$&$24$&$\mathbb{P}(1,3,8,a_3,a_4)$\\
\hline
$E_{12}$&$y^7+z^3$&$(3,7)$&$21$&$42$&$\mathbb{P}(1,6,14,a_3,a_4)$\\
\hline
$W_{13}$&$z(y^4+z^3)$&$(3,4)$&$16$&$16$&$\mathbb{P}(1,3,4,a_3,a_4)$\\
\hline
$Z_{11}$&$y(y^4+z^3)$&$(3,4)$&$15$&$30$&$\mathbb{P}(1,6,8,a_3,a_4)$\\
\hline
$Z_{13}$&$y(y^5+z^3)$&$(3,5)$&$18$&$18$&$\mathbb{P}(1,3,5,a_3,a_4)$\\
\hline
$E_{13}$&$z(y^5+z^2)$&$(2,5)$&$15$&$30$&$\mathbb{P}(1,4,10,a_3,a_4)$\\
\hline
$E_{16}$&$y^9+z^3$&$(1,3)$&$9$&$18$&$\mathbb{P}(1,2,6,a_3,a_4)$\\
\hline
$Y_5$&$y^6+z^6$&$(1,1)$&$6$&$6$&$\mathbb{P}(1,1,1,a_3,a_4)$\\
\hline
$Y_7$&$y^8+z^4$&$(1,2)$&$8$&$8$&$\mathbb{P}(1,1,2,a_3,a_4)$\\
\hline
$Y_{10}$&$y^{12}+z^3$&$(1,4)$&$12$&$12$&$\mathbb{P}(1,1,4,a_3,a_4)$\\
\hline
$Y_{42}$&$y(y^9+z^3)$&$(1,3)$&$10$&$10$&$\mathbb{P}(1,1,3,a_3,a_4)$\\
\hline
\end{tabular}
\end{table}

\begin{remark}
In Lemma--Definition~\ref{def:Z_12-W_15-Z_1,0-N_16-arrangements}, the line $L$ does not pass through the curve singularity $\Sigma$ at $\xi$. On the other hand, locally analytically, we can also obtain the singularities $Z_{12},Z_{15},N_{16}$ using the line $L$ (this does not apply to $W_{15}$). It follows from the proof of Theorem~\ref{thm:classification-of-K3-tails} (Case~1) that, although these alternative singularities produce boundary strata of the KSBA compactification $\overline{\mathrm{M}}_{\mathrm{K}}$, they do not give rise to K3 tails.
\end{remark}




\end{document}